\PassOptionsToPackage{unicode}{hyperref}
\PassOptionsToPackage{hyphens}{url}
\PassOptionsToPackage{dvipsnames,svgnames,x11names}{xcolor}
\documentclass[
]{amsart}
\usepackage{xcolor}
\usepackage{amsmath,amssymb}
\usepackage{iftex}
\ifPDFTeX
  \usepackage[T1]{fontenc}
  \usepackage[utf8]{inputenc}
  \usepackage{textcomp} % provide euro and other symbols
\else % if luatex or xetex
  \usepackage{unicode-math} % this also loads fontspec
  \defaultfontfeatures{Scale=MatchLowercase}
  \defaultfontfeatures[\rmfamily]{Ligatures=TeX,Scale=1}
\fi
\usepackage{lmodern}
\ifPDFTeX\else
\fi
\IfFileExists{upquote.sty}{\usepackage{upquote}}{}
\IfFileExists{microtype.sty}{% use microtype if available
  \usepackage[]{microtype}
  \UseMicrotypeSet[protrusion]{basicmath} % disable protrusion for tt fonts
}{}
\makeatletter
\@ifundefined{KOMAClassName}{% if non-KOMA class
  \IfFileExists{parskip.sty}{%
    \usepackage{parskip}
  }{% else
    \setlength{\parindent}{0pt}
    \setlength{\parskip}{6pt plus 2pt minus 1pt}}
}{% if KOMA class
  \KOMAoptions{parskip=half}}
\makeatother
\makeatletter
\ifx\paragraph\undefined\else
  \let\oldparagraph\paragraph
  \renewcommand{\paragraph}{
    \@ifstar
      \xxxParagraphStar
      \xxxParagraphNoStar
  }
  \newcommand{\xxxParagraphStar}[1]{\oldparagraph*{#1}\mbox{}}
  \newcommand{\xxxParagraphNoStar}[1]{\oldparagraph{#1}\mbox{}}
\fi
\ifx\subparagraph\undefined\else
  \let\oldsubparagraph\subparagraph
  \renewcommand{\subparagraph}{
    \@ifstar
      \xxxSubParagraphStar
      \xxxSubParagraphNoStar
  }
  \newcommand{\xxxSubParagraphStar}[1]{\oldsubparagraph*{#1}\mbox{}}
  \newcommand{\xxxSubParagraphNoStar}[1]{\oldsubparagraph{#1}\mbox{}}
\fi
\makeatother

\usepackage{longtable,booktabs,array}
\usepackage{calc} % for calculating minipage widths
\usepackage{etoolbox}
\makeatletter
\patchcmd\longtable{\par}{\if@noskipsec\mbox{}\fi\par}{}{}
\makeatother
\IfFileExists{footnotehyper.sty}{\usepackage{footnotehyper}}{\usepackage{footnote}}
\makesavenoteenv{longtable}
\usepackage{graphicx}
\makeatletter
\newsavebox\pandoc@box
\newcommand*\pandocbounded[1]{% scales image to fit in text height/width
  \sbox\pandoc@box{#1}%
  \Gscale@div\@tempa{\textheight}{\dimexpr\ht\pandoc@box+\dp\pandoc@box\relax}%
  \Gscale@div\@tempb{\linewidth}{\wd\pandoc@box}%
  \ifdim\@tempb\p@<\@tempa\p@\let\@tempa\@tempb\fi% select the smaller of both
  \ifdim\@tempa\p@<\p@\scalebox{\@tempa}{\usebox\pandoc@box}%
  \else\usebox{\pandoc@box}%
  \fi%
}
\def\fps@figure{htbp}
\makeatother

\providecommand{\tightlist}{%
  \setlength{\itemsep}{0pt}\setlength{\parskip}{0pt}}

\usepackage[numbers,square,sort]{natbib}
\definecolor{linkblue}{RGB}{26,72,140}
\AtBeginDocument{\urlstyle{tt}}
\AtBeginDocument{\address{Ivey Business School, Western University, London, Ontario, Canada}\email{yemine@ivey.ca, bghaddar@ivey.ca}}
\AtBeginDocument{}
\makeatletter\g@addto@macro\appendix{\def\p@section{Appendix~}}\makeatother
\AtBeginDocument{\bibliographystyle{splncs04}}
\AtBeginDocument{\DeclareUrlCommand\doiurl{\urlstyle{rm}}\def\doi{https://doi.org/\doiurl}}
\usepackage{amssymb}
\usepackage{xcolor}
\usepackage{booktabs}
\usepackage{iftex}
\usepackage[section]{placeins}
\usepackage{mleftright}
\usepackage{pgfplots}
\usepackage{siunitx}
\usepackage[l3]{csvsimple}

\pgfplotsset{compat=1.18}
\mleftright

\makeatletter
\AtBeginDocument{%
  \let\sharedthm@master\relax
  \@for\sharedthm@env:=theorem,lemma,proposition,corollary,conjecture,definition,example,exercise,refremark\do{%
    \@ifundefined{c@\sharedthm@env}{}{%
      \ifx\sharedthm@master\relax
        \edef\sharedthm@master{\sharedthm@env}%
      \else
        \expandafter\let\csname c@\sharedthm@env\expandafter\endcsname\csname c@\sharedthm@master\endcsname
        \expandafter\let\csname the\sharedthm@env\expandafter\endcsname\csname the\sharedthm@master\endcsname
      \fi}}}
\makeatother

\AtBeginDocument{\urlstyle{tt}}

\ifLuaTeX
  \AtBeginDocument{\DeclareMathAlphabet{\mathcal}{OMS}{cmsy}{m}{n}}
  \setmathrm{Latin Modern Roman}
  \DeclareFontFamily{U}{msb}{}
  \DeclareFontShape{U}{msb}{m}{n}{<-6>msbm5<6-8>msbm7<8->msbm10}{}
  \AtBeginDocument{\renewcommand{\varnothing}{\mathord{\text{\usefont{U}{msb}{m}{n}\symbol{"3F}}}}}
\fi

\definecolor{notebg}{RGB}{255,246,222}
\definecolor{noteframe}{RGB}{214,158,46}
\definecolor{todofg}{RGB}{160,100,10}
\definecolor{plotblue}{HTML}{2A78D6}
\definecolor{plotorange}{HTML}{EB6834}
\definecolor{plotgreen}{HTML}{1BAF7A}
\definecolor{plotmuted}{HTML}{52514E}

\providecommand{\appendixnumbering}{%
  \renewcommand{\thelemma}{\Alph{section}.\arabic{lemma}}%
  \renewcommand{\theproposition}{\Alph{section}.\arabic{proposition}}%
  \counterwithin*{table}{section}%
  \renewcommand{\thetable}{\Alph{section}.\arabic{table}}}
\providecommand{\coloneqq}{\mathrel{\mathop:}=}

\newcommand{\level}[1]{^{\left(#1\right)}}
\NewDocumentCommand{\DeclareSized}{mm}{\NewDocumentCommand{#1}{o}{#2\IfValueT{##1}{_{##1}}}}
\NewDocumentCommand{\DeclareHierarchy}{mm}{\NewDocumentCommand{#1}{om}{\IfValueTF{##1}{#2[##1]}{#2}\level{##2}}}
\NewDocumentCommand{\DeclareLeveled}{mm}{\NewDocumentCommand{#1}{m}{#2\level{##1}}}
\NewDocumentCommand{\DeclareIndexed}{mm}{\NewDocumentCommand{#1}{m}{#2_{##1}}}

\newcommand{\Reals}{\mathbb{R}}
\newcommand{\Naturals}{\mathbb{N}}

\newcommand{\GenericCone}{\mathcal{K}}
\DeclareSized{\COP}{\mathcal{C}^{*}}
\DeclareSized{\CP}{\mathcal{C}}
\DeclareSized{\DNN}{\mathcal{D}}
\DeclareSized{\PSD}{\mathcal{S}^{+}}
\DeclareSized{\Sym}{\mathcal{S}}
\DeclareSized{\NN}{\mathcal{N}}
\DeclareSized{\POL}{\mathcal{P}}
\DeclareSized{\SOS}{\mathcal{R}}
\newcommand{\widecone}[1]{\widetilde{#1\mkern2mu}\mkern-2mu}
\DeclareSized{\DSOS}{\widecone{\POL}}
\DeclareSized{\SDSOS}{\widecone{\SOS}}

\DeclareHierarchy{\PolyaCone}{\POL}
\DeclareHierarchy{\PolyaSOSCone}{\SOS}
\DeclareHierarchy{\PolyaDSOSCone}{\DSOS}
\DeclareHierarchy{\PolyaSDSOSCone}{\SDSOS}

\newcommand{\FormVariable}{x}
\newcommand{\forms}[1]{#1\left[\FormVariable\right]}
\newcommand{\NonnegativeForms}{\forms{\Reals_{+}}}
\newcommand{\DSOSForms}{\forms{\Sigma_{\pm}}}
\newcommand{\SDSOSForms}{\forms{\Sigma_{2}}}
\newcommand{\SOSForms}{\forms{\Sigma}}

\newcommand{\Biclique}{g}
\newcommand{\BalancedBiclique}{\Biclique_{\mathrm{bal}}}
\newcommand{\BicliqueRatio}{h}
\newcommand{\Independence}{\alpha}
\newcommand{\Biindependence}{\Independence^{*}}
\newcommand{\LovaszTheta}{\vartheta}
\newcommand{\SchrijverTheta}{\LovaszTheta^{\prime}}
\newcommand{\Complement}[1]{\widehat{#1}}

\DeclareSized{\Identity}{I}
\DeclareSized{\OnesMatrix}{J}
\DeclareSized{\OnesVector}{\mathbf{1}}
\newcommand{\Biadjacency}{U}
\newcommand{\Adjacency}{W}
\newcommand{\CostMatrix}{C}

\DeclareIndexed{\LasserreBicliqueRatio}{\BicliqueRatio}
\DeclareIndexed{\LasserreBiclique}{\Biclique}
\newcommand{\ConeBicliqueRatio}[1]{\LasserreBicliqueRatio{1}^{#1}}
\newcommand{\SudaTanakaBound}{T}
\DeclareLeveled{\Polya}{\zeta}
\DeclareLeveled{\PolyaDSOS}{\tilde{\zeta}}
\DeclareLeveled{\PolyaSOS}{\rho}
\DeclareLeveled{\PolyaSDSOS}{\tilde{\rho}}

\newcommand{\EdgeMultiplier}{w}

\DeclareIndexed{\PenaltyMatrix}{H}
\newcommand{\BilinearValue}{\sigma}
\DeclareIndexed{\BilinearMax}{\BilinearValue}
\DeclareIndexed{\QuadraticMax}{\mu}
\newcommand{\SmallestPenalty}{\gamma^{*}}
\DeclareIndexed{\PenaltyThreshold}{\gamma}

\newcommand{\Transpose}{T}

\DeclareMathOperator{\tr}{\mathbf{tr}}
\DeclareMathOperator{\nz}{\mathbf{nz}}
\DeclareMathOperator{\Diag}{Diag}
\DeclareMathOperator{\diag}{diag}
\DeclareMathOperator{\sgn}{sgn}

\newcommand{\TableData}[1]{tables/#1.csv}
\newcommand{\FigureData}[1]{figures/#1.csv}

\csvreader[head to column names=false]{\TableData{values}}{name=\valuename,value=\valuevalue}{%
  \expandafter\xdef\csname value@\valuename\endcsname{\valuevalue}}
\newcommand{\Value}[1]{%
  \ifcsname value@#1\endcsname
    \csname value@#1\endcsname
  \else
    \errmessage{no value #1 in values.csv}%
  \fi}

\newcommand{\allof}[2]{\ifnum#1=#2 all $#2$\else $#1$ of the $#2$\fi}

\NewDocumentCommand{\fixed}{mm}{\num[round-mode=places,round-precision=#1,round-half=up]{#2}}
\NewDocumentCommand{\fixedtrim}{mm}{\num{\fpeval{round(#2,#1)}}}
\NewDocumentCommand{\atmost}{mm}{\fixed{#1}{\fpeval{ceil(#2,#1)}}}
\NewDocumentCommand{\atleast}{mm}{\fixed{#1}{\fpeval{floor(#2,#1)}}}

\NewDocumentCommand{\sci}{m}{\num[print-unity-mantissa=false,exponent-product=\cdot]{#1}}

\ExplSyntaxOn
\fp_new:N \l__paper_exponent_fp
\fp_new:N \l__paper_mantissa_fp
\cs_new_protected:Npn \__paper_sci:nnn #1#2#3
  {
    \fp_set:Nn \l__paper_exponent_fp { floor(ln(#2) / ln(10)) }
    \fp_set:Nn \l__paper_mantissa_fp { #3(#2 / 10^(\l__paper_exponent_fp), #1 - 1) }
    \fp_compare:nNnT { \l__paper_mantissa_fp } = { 10 }
      {
        \fp_set:Nn \l__paper_mantissa_fp { 1 }
        \fp_add:Nn \l__paper_exponent_fp { 1 }
      }
    \sci { \fp_use:N \l__paper_mantissa_fp e \fp_use:N \l__paper_exponent_fp }
  }
\NewDocumentCommand{\sciatmost}{mm}{ \__paper_sci:nnn {#1} {#2} {ceil} }
\NewDocumentCommand{\sciatleast}{mm}{ \__paper_sci:nnn {#1} {#2} {floor} }
\ExplSyntaxOff

\NewDocumentCommand{\readrows}{mmmm}{%
  \csvreader[head to column names, head to column names prefix=col, collect data,
    late after line=\\, late after last line=\\, #2]{#3}{}{#4}%
  \let#1\csvdatacollection}

\ExplSyntaxOn
\NewDocumentCommand{\gapcell}{mm}
  {
    \fp_compare:nNnTF {#2} = {0}
      { \rlap{$0$} \phantom{\fixed{#1}{0}} }
      {
        \fp_compare:nNnTF {#2} < {10^(-#1)}
          { {<} \fixed{#1}{\fp_eval:n {10^(-#1)}} }
          { \fixed{#1}{#2} }
      }
  }
\ExplSyntaxOff

\newcommand{\boundsheads}{%
  \# & $\LasserreBicliqueRatio{1}$ & $\ConeBicliqueRatio{\DNN}$ & $\PolyaSOS{1}$
  & $\LasserreBicliqueRatio{1}$ & $\ConeBicliqueRatio{\DNN}$ & $\PolyaSOS{1}$
  & $\LasserreBicliqueRatio{1}$ & $\ConeBicliqueRatio{\DNN}$ & $\PolyaSOS{1}$}
\newcommand{\boundscells}{%
  $\colcount$ & $\colexacthone$ & $\colexacthonednn$ & $\colexactrhoone$
  & $\gapcell{1}{\colmeanhone}$ & $\gapcell{2}{\colmeanhonednn}$ & $\gapcell{3}{\colmeanrhoone}$
  & $\gapcell{1}{\colmaxhone}$ & $\gapcell{1}{\colmaxhonednn}$ & $\gapcell{3}{\colmaxrhoone}$}
\makeatletter
\@ifpackageloaded{caption}{}{\usepackage{caption}}
\AtBeginDocument{%
\ifdefined\contentsname
  \renewcommand*\contentsname{Table of contents}
\else
  \newcommand\contentsname{Table of contents}
\fi
\ifdefined\listfigurename
  \renewcommand*\listfigurename{List of Figures}
\else
  \newcommand\listfigurename{List of Figures}
\fi
\ifdefined\listtablename
  \renewcommand*\listtablename{List of Tables}
\else
  \newcommand\listtablename{List of Tables}
\fi
\ifdefined\figurename
  \renewcommand*\figurename{Figure}
\else
  \newcommand\figurename{Figure}
\fi
\ifdefined\tablename
  \renewcommand*\tablename{Table}
\else
  \newcommand\tablename{Table}
\fi
}
\@ifpackageloaded{float}{}{\usepackage{float}}
\floatstyle{ruled}
\@ifundefined{c@chapter}{\newfloat{codelisting}{h}{lop}}{\newfloat{codelisting}{h}{lop}[chapter]}
\floatname{codelisting}{Listing}

\usepackage{amsthm}
\theoremstyle{plain}
\newtheorem{corollary}{Corollary}[section]
\theoremstyle{plain}
\newtheorem{lemma}{Lemma}[section]
\theoremstyle{plain}
\newtheorem{proposition}{Proposition}[section]
\theoremstyle{plain}
\newtheorem{theorem}{Theorem}[section]
\theoremstyle{remark}
\AtBeginDocument{}

\newtheorem{refremark}{Remark}[section]

\makeatother
\makeatletter
\@ifpackageloaded{caption}{}{\usepackage{caption}}
\@ifpackageloaded{subcaption}{}{\usepackage{subcaption}}
\makeatother
\usepackage{bookmark}
\IfFileExists{xurl.sty}{\usepackage{xurl}}{} % add URL line breaks if available
\makeatletter
\@ifundefined{xmpquote}{\newcommand{\xmpquote}[1]{#1}}{}
\makeatother
\hypersetup{
  pdftitle={Exact Penalties for Maximum-Edge Bicliques in Bipartite Graphs},
  pdfauthor={Youssouf Emine; Bissan Ghaddar},
  pdfkeywords={\xmpquote{maximum-edge biclique}, \xmpquote{penalty formulation}, \xmpquote{biindependent pair}, \xmpquote{copositive and completely positive matrices}, \xmpquote{Motzkin--Straus}},
  colorlinks=true,
  linkcolor={linkblue},
  filecolor={Maroon},
  citecolor={linkblue},
  urlcolor={linkblue},
  pdfcreator={LaTeX via pandoc}}

\title{Exact Penalties for Maximum-Edge Bicliques in Bipartite Graphs}
\author{Youssouf Emine \and Bissan Ghaddar}
\date{}
\begin{document}
\begin{abstract}
The maximum-edge biclique problem asks for a complete bipartite subgraph with the most edges in a bipartite graph, and it is NP-hard.
In the bipartite complement, it asks for two vertex sets, one on each side and with no edge between them, whose sizes have the largest product \(g\).
A known nonconvex quadratic program penalizes the edges between the two sets.
Its global maximizers are the optimal pairs for large penalties, and the penalties known to suffice grow with the size of the graph.
We show that every penalty of at least \(\frac12\left(1+\sqrt{g}\right)\) suffices, and that the local maximizers with positive value are then pairs that cannot be enlarged.
This threshold is sharp up to an additive one half, and it is never larger than the earlier thresholds.
The same program is a copositive formulation that gives the optimal value exactly, and so does its completely positive dual.
Tractable relaxations of these programs bound the unknown \(g\) and give computable penalties.
The doubly nonnegative relaxation is inexact already on six vertices, and Pólya's linear programming bounds are never exact.
On random graphs with 30 vertices and edge probability 0.3, a local ascent from random starting points reached an optimal pair on \(\Value{ascentgraphs}\) of 20 graphs with the threshold penalty and on none with the two earlier thresholds.
\end{abstract}
\maketitle

\section{Introduction}\label{sec-intro}

The maximum-edge biclique problem asks for a complete bipartite subgraph, or biclique, with the largest number of edges in a bipartite graph, and it is NP-hard \citep{peeters2003maximum}.
Equivalently, it asks for a largest all-ones submatrix of a binary matrix, a task that arises in the biclustering of gene expression data \citep{prelic2006systematic} and in data mining \citep{geerts2004tiling, liu2006efficient} and is closely related to rank-one nonnegative factorization \citep{gillis2008nonnegative}.
Let \(G=\left(V_{1}\cup V_{2},E\right)\) be a bipartite graph with vertex set \(V=V_{1}\cup V_{2}\), \(m_{1}=|V_{1}|\ge1\), \(m_{2}=|V_{2}|\ge1\) and \(n=m_{1}+m_{2}\).
A pair \(\left(A,B\right)\) with \(A\subseteq V_{1}\) and \(B\subseteq V_{2}\) is \emph{biindependent} if no edge joins \(A\) and \(B\), and
\[
\Biclique\left(G\right)\coloneqq\max\left\{|A||B|:\left(A,B\right)\ \text{biindependent}\right\}.
\]
The biindependent pairs of \(G\) are the bicliques of its bipartite complement \(\Complement{G}\), which has the same two sides and, between them, exactly the edges that \(G\) lacks.
Computing \(\Biclique\left(G\right)\) is therefore the maximum-edge biclique problem on \(\Complement{G}\).
A biindependent pair is \emph{maximal} if no vertex can be added to either side and \emph{maximum} if \(|A||B|=\Biclique\left(G\right)\), and we say maximal pair and maximum pair for short.

For the maximum clique problem, the Motzkin--Straus theorem \citep{motzkin1965maxima} turns a combinatorial problem into the maximization of a quadratic form over the simplex, and regularized or penalized versions of it make the local maximizers correspond to maximal cliques \citep{hungerford2019general, belachew2017solving}.
Continuous quadratic formulations of this kind are known for several other graph problems \citep{hager2015continuous}.
For the maximum-edge biclique problem, Gillis and Glineur \citep{gillis2014continuous} maximize the bilinear form \(y^{\Transpose}\left(\OnesMatrix-2\gamma \Biadjacency\right)z\) over nonnegative unit vectors \(y\) and \(z\).
Here \(\Biadjacency\in\left\{0,1\right\}^{V_{1}\times V_{2}}\) is the biadjacency matrix of \(G\) and \(\OnesMatrix\in\Reals^{V_{1}\times V_{2}}\) the all-ones matrix, and \(\gamma\ge0\) penalizes the edges of \(G\) between the supports of \(y\) and \(z\).
Gillis and Glineur write this formulation for the graph \(\Complement{G}\), whose bicliques are our biindependent pairs, and with the parameter \(2\gamma-1\) where we use \(\gamma\).
We quote the earlier results in our notation.

We study this formulation in a symmetric form.
With rows and columns indexed by \(V_{1}\) and then \(V_{2}\), let
\[
\CostMatrix\coloneqq\frac12\begin{pmatrix}0&\OnesMatrix\\ \OnesMatrix^{\Transpose}&0\end{pmatrix}\in\Sym[n],\qquad\text{so that}\qquad u^{\Transpose}\CostMatrix u=\left(\sum_{i\in V_{1}}u_{i}\right)\left(\sum_{j\in V_{2}}u_{j}\right).
\]
Let \(\Adjacency\) be the adjacency matrix of \(G\).
For \(\gamma\ge0\), the nonconvex quadratic program
\[
\max\left\{u^{\Transpose}\left(\CostMatrix-\gamma \Adjacency\right)u:u\in\Reals^{V}_{+},\ \|u\|=1\right\}
\]
is equivalent to the bilinear formulation (\hyperref[sec-forms]{Section~\ref*{sec-forms}}).
We call the penalty \(\gamma\) \emph{exact} if the maximum of this program is \(\tfrac12\sqrt{\Biclique\left(G\right)}\), the value of a maximum pair.

The penalties known to be exact grow with the size of the graph.
Gillis and Glineur prove that the bilinear form has the maximum pairs as global maximizers and the maximal pairs as local maximizers once \(2\gamma-1\ge\max\left(m_{1},m_{2}\right)\) \citep[Theorems 1 and 2]{gillis2014continuous}, and that \(2\gamma-1\ge\sqrt{\left|E\left(\Complement{G}\right)\right|}\) already makes every global maximizer a maximum pair \citep[Lemma 3]{gillis2014continuous}.
For cliques, the counterpart of this formulation due to Belachew and Gillis \citep[Theorems 3 and 4]{belachew2017solving} requires a penalty of at least the number of vertices.
A large penalty hurts the local search methods that work directly on these formulations \citep{gillis2014continuous, belachew2017solving}, because each ascent step then discards more vertices at once and the search stops at a maximal pair near a small support (\hyperref[sec-ascent]{Section~\ref*{sec-ascent}}).
We ask how small an exact penalty can be.

Our main result is that the exact penalties are set by the optimal value, not by the size of the graph.
Every \(\gamma\ge\frac12\left(1+\sqrt{\Biclique\left(G\right)}\right)\) is exact, and if \(\Biclique\left(G\right)\ge1\) the maximizers of the quadratic program are then exactly the vectors that are uniform on each side of a maximum pair.
The reason is that a local maximizer of the bilinear form with positive value and an edge of \(G\) between its supports has value larger than \(2\gamma-1\), whereas every point without such an edge has value at most \(\sqrt{\Biclique\left(G\right)}\).
For the same \(\gamma\), the local maximizers of the bilinear form with positive value are uniform on maximal pairs.
The threshold cannot be lowered by more than \(\tfrac12\) in general, and in the parameter \(2\gamma-1\) it reads \(\sqrt{\Biclique\left(G\right)}\), which is never larger than the two thresholds of Gillis and Glineur.
It involves the unknown \(\Biclique\left(G\right)\), but any upper bound can take its place, and we obtain upper bounds from the largest singular value of the biadjacency matrix of \(\Complement{G}\) and from semidefinite programs.

The semidefinite programs come from a conic view of the penalized program.
With an exact penalty, this program is an exact copositive formulation of \(\tfrac12\sqrt{\Biclique\left(G\right)}\) in which every edge of \(G\) carries the same multiplier \(\gamma\).
With free multipliers, it is the dual of a completely positive program that strengthens a semidefinite bound of Laurent et al. \citep{laurent2025semidefinite}.
We prove that the completely positive program is exact, and our threshold then shows that a uniform multiplier attains the dual optimum, which conic duality alone does not guarantee.
Tractable relaxations of the two exact programs give upper bounds on \(\tfrac12\sqrt{\Biclique\left(G\right)}\), and hence computable exact penalties.

The experiments use graphs with at most \(40\) vertices.
On the random graphs of \hyperref[sec-bounds]{Section~\ref*{sec-bounds}}, a rounded semidefinite bound already equals \(\Biclique\left(G\right)\), so that the penalty it gives is the threshold itself and serves to find a maximum pair, not to compute the optimal value.
From random starting points on the \(20\) graphs with \(15+15\) vertices and edge probability \(0.3\), a local ascent on the bilinear form reached a maximum pair on \(\Value{ascentgraphs}\) graphs with the threshold and on none with the three larger thresholds of \hyperref[sec-exp-threshold]{Section~\ref*{sec-exp-threshold}}.

Our contributions are the following.

\begin{enumerate}
\def\labelenumi{\arabic{enumi}.}
\tightlist
\item
  The exact penalty threshold (\hyperref[thm-penalty]{Theorem~\ref*{thm-penalty}}), its extension to the local maximizers (\hyperref[prp-local]{Proposition~\ref*{prp-local}}) and its sharpness (\hyperref[prp-sharp]{Proposition~\ref*{prp-sharp}}).
\item
  Exact completely positive and copositive formulations of \(\tfrac12\sqrt{\Biclique\left(G\right)}\) (\hyperref[thm-cp]{Theorem~\ref*{thm-cp}}, \hyperref[cor-dual]{Corollary~\ref*{cor-dual}}).
\item
  Three relaxations: the doubly nonnegative bound is exact up to five vertices but not on six (\hyperref[sec-dnn]{Section~\ref*{sec-dnn}}), the level-1 semidefinite bound of Parrilo is never weaker (\hyperref[sec-rho]{Section~\ref*{sec-rho}}), and the Pólya linear programming bounds are finite from an explicit level and never exact (\hyperref[sec-polya]{Section~\ref*{sec-polya}}).
\item
  Computable penalties and their effect on a local ascent (\hyperref[sec-computable]{Section~\ref*{sec-computable}}, \hyperref[sec-numerics]{Section~\ref*{sec-numerics}}).
\end{enumerate}

\textbf{Relation to earlier work.}
The conic results follow the approach of de Klerk and Pasechnik \citep{deklerk2002approximation} for the stability number, refined by Peña et al. \citep{pena2007computing}, and the differences are described in \hyperref[sec-cp]{Section~\ref*{sec-cp}} and \hyperref[sec-polya]{Section~\ref*{sec-polya}}.
We found no earlier copositive formulation of \(\Biclique\left(G\right)\) or \(\tfrac12\sqrt{\Biclique\left(G\right)}\), apart from the lifted one in dimension \(1+2n\) that the theorem of Burer \citep{burer2009copositive} gives for every binary quadratic program.

\textbf{Paper outline.}
In \hyperref[sec-prelim]{Section~\ref*{sec-prelim}} we fix the notation and recall the semidefinite bounds of Laurent et al.
We prove the threshold in \hyperref[sec-penalty]{Section~\ref*{sec-penalty}}, derive the exact conic formulations in \hyperref[sec-cp]{Section~\ref*{sec-cp}} and study their relaxations in \hyperref[sec-relax]{Section~\ref*{sec-relax}}.
The experiments are in \hyperref[sec-numerics]{Section~\ref*{sec-numerics}}, and \hyperref[sec-conclusion]{Section~\ref*{sec-conclusion}} concludes.

\section{Preliminaries}\label{sec-prelim}

\textbf{Graphs and vectors.}
We identify \(V\) with \(\left\{1,\dots,n\right\}\), listing \(V_{1}\) first, so that the entries of vectors in \(\Reals^{n}\) and the rows and columns of \(n\times n\) matrices are indexed by \(V\).
We write the size of \(G\) as the sum \(m_{1}+m_{2}\), so that a graph with \(3+4\) vertices has three vertices in \(V_{1}\) and four in \(V_{2}\).
A set \(S\subseteq V\) is stable in \(G\) if and only if \(\left(S\cap V_{1},S\cap V_{2}\right)\) is biindependent.
If \(G\) is complete bipartite, then \(\Biclique\left(G\right)=0\), since one of \(A\), \(B\) must be empty, and otherwise \(\Biclique\left(G\right)\ge1\).
We omit the argument \(G\) of a parameter or a bound when the graph is clear from the context.

The maximizers of the penalized program are written with the unit vectors of sets.
For \(S\subseteq V\), \(\OnesVector[S]\) is its indicator vector, and if \(S\ne\varnothing\), \(e_{S}\coloneqq\frac1{\sqrt{|S|}}\OnesVector[S]\) is its unit vector, so that \(e_{\left\{i\right\}}=e_{i}\).
We use \(\OnesVector[S]\) and \(e_{S}\) also as vectors of \(\Reals^{V_{1}}\) or \(\Reals^{V_{2}}\) when \(S\) lies in one side.
For a vector \(u\), \(\nz\left(u\right)\coloneqq\left\{i:u_{i}\ne0\right\}\) is its support and \(\Diag\left(u\right)\) is the diagonal matrix with diagonal \(u\), while for a square matrix \(M\), \(\diag M\) is the vector of its diagonal entries.

\textbf{Matrices and cones.}
We write \(\Sym[n]\) for the space of real symmetric \(n\times n\) matrices, with the inner product \(\left\langle A,B\right\rangle=\tr\left(AB\right)\), and \(\PSD[n]\) and \(\NN[n]\) for its cones of positive semidefinite and of entrywise nonnegative matrices.
The conic programs of \hyperref[sec-cp]{Section~\ref*{sec-cp}} and \hyperref[sec-relax]{Section~\ref*{sec-relax}} use the doubly nonnegative, completely positive and copositive cones
\[
\DNN[n]\coloneqq\PSD[n]\cap\NN[n],\qquad \CP[n]\coloneqq\operatorname{conv}\left\{uu^{\Transpose}:u\in\Reals^{n}_{+}\right\},
\]
and \(\COP[n]\coloneqq\left\{M\in\Sym[n]:u^{\Transpose}Mu\ge0\ \forall u\in\Reals^{n}_{+}\right\}\).
The notation \(\COP[n]\) reflects that the copositive cone is the dual of \(\CP[n]\), and we have \(\CP[n]\subseteq\DNN[n]\subseteq\PSD[n]\).
The matrix \(\OnesMatrix-\Biadjacency\) is the biadjacency matrix of \(\Complement{G}\).
Writing \(u=\left(y,z\right)\) with \(y\in\Reals^{V_{1}}\) and \(z\in\Reals^{V_{2}}\), the adjacency matrix and the two quadratic forms are
\begin{equation}\protect\phantomsection\label{eq-split}{
\Adjacency\coloneqq\begin{pmatrix}0&\Biadjacency\\\Biadjacency^{\Transpose}&0\end{pmatrix},\qquad u^{\Transpose}\CostMatrix u=y^{\Transpose}\OnesMatrix z,\qquad u^{\Transpose}\Adjacency u=2y^{\Transpose}\Biadjacency z .
}\end{equation}

\textbf{Bounds of Laurent et al.}
The semidefinite bounds that give computable penalties in \hyperref[sec-computable]{Section~\ref*{sec-computable}} come from Laurent et al. \citep{laurent2025semidefinite}.
Besides \(\Biclique\), that paper studies the largest ratio \(\BicliqueRatio\left(G\right)\) of \(|A||B|\) to \(|A|+|B|\) over biindependent pairs with \(A\cup B\ne\varnothing\), which is NP-hard as well and satisfies \(\BicliqueRatio\left(G\right)\le\tfrac12\sqrt{\Biclique\left(G\right)}\) \citep[Lemma 1]{laurent2025semidefinite}.
The first of the semidefinite bounds of that paper is
\begin{equation}\protect\phantomsection\label{eq-h1}{
\LasserreBicliqueRatio{1}\left(G\right)\coloneqq\max\left\{\left\langle \CostMatrix,X\right\rangle:X\in\PSD[n],\ \tr\left(X\right)=1,\ X_{ij}=0\ \left(ij\in E\right)\right\},
}\end{equation}
and we write \(X_{E}=0\) for the last constraint.
It satisfies \(\tfrac12\sqrt{\Biclique\left(G\right)}\le \LasserreBicliqueRatio{1}\left(G\right)\) \citep[Proposition 1]{laurent2025semidefinite} and is the analogue of the theta number for biindependent pairs.
To compare cones, we write, for a closed convex cone \(\GenericCone\subseteq\Sym[n]\),
\begin{equation}\protect\phantomsection\label{eq-h1K}{
\ConeBicliqueRatio{\GenericCone}\left(G\right)\coloneqq\max\left\{\left\langle \CostMatrix,X\right\rangle:X\in\GenericCone,\ \tr\left(X\right)=1,\ X_{E}=0\right\},
}\end{equation}
which is \(\LasserreBicliqueRatio{1}\left(G\right)\) for \(\GenericCone=\PSD[n]\).

\section{The Penalty Formulation}\label{sec-penalty}

\subsection{Two Forms of the Program}\label{sec-forms}

The penalized program has a symmetric form, which is a copositive program, and a bilinear form, on which the proofs and the local search work.
For \(\gamma\ge0\) let
\begin{equation}\protect\phantomsection\label{eq-mu}{
\QuadraticMax{\gamma}\left(G\right)\coloneqq\max\left\{u^{\Transpose}\left(\CostMatrix-\gamma \Adjacency\right)u:u\in\Reals^{V}_{+},\ \|u\|=1\right\}.
}\end{equation}
Since, for a symmetric matrix \(Q\), \(\lambda \Identity-Q\in\COP[n]\) if and only if \(u^{\Transpose}Qu\le\lambda\|u\|^{2}\) for all \(u\ge0\),
\begin{equation}\protect\phantomsection\label{eq-mu-cop}{
\QuadraticMax{\gamma}\left(G\right)=\min\left\{\lambda:\lambda \Identity-\CostMatrix+\gamma \Adjacency\in\COP[n]\right\}.
}\end{equation}
By \hyperref[eq-split]{(\ref*{eq-split})}, \(u^{\Transpose}\left(\CostMatrix-\gamma \Adjacency\right)u=y^{\Transpose}\PenaltyMatrix{\gamma} z\) with \(\PenaltyMatrix{\gamma}\coloneqq\OnesMatrix-2\gamma \Biadjacency\).
Normalizing \(y\) and \(z\) separately, we put
\begin{equation}\protect\phantomsection\label{eq-sigma}{
\BilinearMax{\gamma}\left(G\right)\coloneqq\max\left\{y^{\Transpose}\PenaltyMatrix{\gamma} z:y\in\Reals^{V_{1}}_{+},\ z\in\Reals^{V_{2}}_{+},\ \|y\|=\|z\|=1\right\}.
}\end{equation}
Writing a feasible \(u\) of \hyperref[eq-mu]{(\ref*{eq-mu})} as \(\left(a\hat y,b\hat z\right)\) with \(a^{2}+b^{2}=1\) and unit vectors \(\hat y\), \(\hat z\) (arbitrary if \(y=0\) or \(z=0\)), we get \(u^{\Transpose}\left(\CostMatrix-\gamma \Adjacency\right)u=ab\,\hat y^{\Transpose}\PenaltyMatrix{\gamma}\hat z\) with \(ab\in\left[0,\frac12\right]\), and hence
\begin{equation}\protect\phantomsection\label{eq-mu-sigma}{
\QuadraticMax{\gamma}\left(G\right)=\tfrac12\max\left\{\BilinearMax{\gamma}\left(G\right),0\right\}.
}\end{equation}
The matrix \(\PenaltyMatrix{\gamma}\) has entries \(1\) on the edges of \(\Complement{G}\) and \(1-2\gamma\) on the edges of \(G\).
Program \hyperref[eq-sigma]{(\ref*{eq-sigma})} is therefore the program (25) of Gillis and Glineur \citep{gillis2014continuous} for \(\Complement{G}\), with parameter \(2\gamma-1\).

\subsection{The Threshold}\label{sec-threshold}

The threshold comes from comparing two kinds of feasible points of \hyperref[eq-sigma]{(\ref*{eq-sigma})}.
At a local maximizer whose supports are joined by an edge of \(G\), moving weight between the two ends of that edge cannot increase the value, and this second-order condition forces a positive value to exceed \(2\gamma-1\) (\hyperref[lem-second-order]{Lemma~\ref*{lem-second-order}}).
A point whose supports are joined by no edge has value at most \(\sqrt{\Biclique}\) by the Cauchy--Schwarz inequality (\hyperref[lem-sigma]{Lemma~\ref*{lem-sigma}}).
As long as the maximum exceeds \(\sqrt{\Biclique}\), a maximizer is therefore of the first kind and the maximum also exceeds \(2\gamma-1\).
Since the maximum is continuous and nonincreasing in \(\gamma\), it must have come down to \(\sqrt{\Biclique}\) by the time \(2\gamma-1\) reaches \(\sqrt{\Biclique}\).
The first step moves weight as the arguments of Hungerford and Rinaldi \citep[Section 2]{hungerford2019general} for regularized Motzkin--Straus formulations do between the two ends of a non-edge.

\begin{lemma}[]\protect\hypertarget{lem-second-order}{}\label{lem-second-order}

Let \(\gamma\ge0\), let \(\left(y,z\right)\) be a local maximizer of \hyperref[eq-sigma]{(\ref*{eq-sigma})} with value \(\BilinearValue=y^{\Transpose}\PenaltyMatrix{\gamma} z\), and let \(S=\nz\left(y\right)\) and \(R=\nz\left(z\right)\).
Then \(\left(\PenaltyMatrix{\gamma} z\right)_{S}=\BilinearValue y_{S}\) and \(\left(\PenaltyMatrix{\gamma}^{\Transpose}y\right)_{R}=\BilinearValue z_{R}\).
If \(i\in S\), \(j\in R\) and \(ij\in E\), then
\[
2\gamma-1\le\BilinearValue\left(1-\tfrac12\left(y_{i}+z_{j}\right)^{2}\right),
\]
and in particular \(2\gamma-1<\BilinearValue\) if \(\BilinearValue>0\).

\end{lemma}

\begin{proof}
For \(\varphi\in\Reals^{V_{1}}\) supported in \(S\) with \(\varphi^{\Transpose}y=0\) and \(\psi\in\Reals^{V_{2}}\) supported in \(R\) with \(\psi^{\Transpose}z=0\), put
\[
y\left(t\right)=\frac{y+t\varphi}{\|y+t\varphi\|},\qquad z\left(t\right)=\frac{z+t\psi}{\|z+t\psi\|}.
\]
Here \(\varphi\) or \(\psi\) may be \(0\), as it must be when \(|S|=1\) or \(|R|=1\).
For \(|t|\) small these vectors are nonnegative, since \(y\) is positive on \(S\), \(z\) is positive on \(R\), and the other entries stay \(0\); hence they are feasible for \hyperref[eq-sigma]{(\ref*{eq-sigma})}, and \(f\left(t\right)=y\left(t\right)^{\Transpose}\PenaltyMatrix{\gamma} z\left(t\right)\) has a local maximum at \(t=0\).
As \(\|y+t\varphi\|^{2}=1+t^{2}\|\varphi\|^{2}\) and \(\|z+t\psi\|^{2}=1+t^{2}\|\psi\|^{2}\), expanding gives
\[
f\left(t\right)=\BilinearValue+t\left(\varphi^{\Transpose}\PenaltyMatrix{\gamma} z+y^{\Transpose}\PenaltyMatrix{\gamma}\psi\right)+t^{2}\left(\varphi^{\Transpose}\PenaltyMatrix{\gamma}\psi-\tfrac\BilinearValue2\left(\|\varphi\|^{2}+\|\psi\|^{2}\right)\right)+O\left(t^{3}\right),
\]
and thus the coefficient of \(t\) vanishes, which is the first-order condition, and that of \(t^{2}\) is nonpositive, which is the second-order condition.
With \(\psi=0\), the first-order condition says that \(\varphi^{\Transpose}\PenaltyMatrix{\gamma} z=0\) for every \(\varphi\) supported in \(S\) with \(\varphi^{\Transpose}y=0\), and hence \(\left(\PenaltyMatrix{\gamma} z\right)_{S}\) is a multiple of \(y_{S}\), and multiplying by \(y\) shows that \(\left(\PenaltyMatrix{\gamma} z\right)_{S}=\BilinearValue y_{S}\); likewise \(\left(\PenaltyMatrix{\gamma}^{\Transpose}y\right)_{R}=\BilinearValue z_{R}\).
Now take \(\varphi=e_{i}-y_{i}y\) and \(\psi=-\left(e_{j}-z_{j}z\right)\).
Then \(\|\varphi\|^{2}=1-y_{i}^{2}\) and \(\|\psi\|^{2}=1-z_{j}^{2}\), and the first-order conditions and \(\left(\PenaltyMatrix{\gamma}\right)_{ij}=1-2\gamma\) give
\[
\varphi^{\Transpose}\PenaltyMatrix{\gamma}\psi=-\left(\PenaltyMatrix{\gamma}\right)_{ij}+\BilinearValue y_{i}z_{j}=2\gamma-1+\BilinearValue y_{i}z_{j} .
\]
The second-order condition becomes \(2\gamma-1+\BilinearValue y_{i}z_{j}\le\frac\BilinearValue2\left(2-y_{i}^{2}-z_{j}^{2}\right)\), that is, \(2\gamma-1\le\BilinearValue\left(1-\frac12\left(y_{i}+z_{j}\right)^{2}\right)\).
If \(\BilinearValue>0\), the right side is less than \(\BilinearValue\) because \(\left(y_{i}+z_{j}\right)^{2}>0\), as \(y_{i}\) and \(z_{j}\) are positive.
\end{proof}

Besides the bound for the points without such an edge, \hyperref[lem-sigma]{Lemma~\ref*{lem-sigma}} gives the convexity and the monotonicity of \(\gamma\mapsto\BilinearMax{\gamma}\) that the continuity argument needs.

\begin{lemma}[]\protect\hypertarget{lem-sigma}{}\label{lem-sigma}

Let \(\Biclique=\Biclique\left(G\right)\ge1\).
The function \(\gamma\mapsto\BilinearMax{\gamma}\left(G\right)\) is convex, continuous and nonincreasing on \(\left[0,\infty\right)\), and \(\BilinearMax{\gamma}\ge\sqrt \Biclique\) for all \(\gamma\).
For every feasible \(\left(y,z\right)\) such that no edge joins \(\nz\left(y\right)\) and \(\nz\left(z\right)\), we have \(y^{\Transpose}\PenaltyMatrix{\gamma} z\le\sqrt \Biclique\), with equality if and only if \(y\) and \(z\) are uniform on their supports and these form a biindependent pair with product \(\Biclique\).

\end{lemma}

\begin{proof}
Fix a feasible \(\left(y,z\right)\).
The function \(\gamma\mapsto y^{\Transpose}\PenaltyMatrix{\gamma} z=y^{\Transpose}\OnesMatrix z-2\gamma\,y^{\Transpose}\Biadjacency z\) is affine, and its slope \(-2y^{\Transpose}\Biadjacency z\) lies in \(\left[-2\|\Biadjacency\|_{2},0\right]\).
Thus \(\gamma\mapsto\BilinearMax{\gamma}\), the maximum of these functions over the compact feasible set, is convex and nonincreasing, and continuous, since it is Lipschitz with constant \(2\|\Biadjacency\|_{2}\).

For the bound \(\BilinearMax{\gamma}\ge\sqrt \Biclique\), take a maximum pair \(\left(A,B\right)\), and let \(y=e_{A}\) and \(z=e_{B}\).
No edge joins \(A\) and \(B\); hence \(y^{\Transpose}\Biadjacency z=0\), and \(y^{\Transpose}\PenaltyMatrix{\gamma} z=y^{\Transpose}\OnesMatrix z=\sqrt{|A||B|}=\sqrt \Biclique\).

For the last statement, let \(S=\nz\left(y\right)\) and \(R=\nz\left(z\right)\).
If no edge joins them, then \(y^{\Transpose}\Biadjacency z=0\) and \(\left(S,R\right)\) is biindependent; thus \(|S||R|\le \Biclique\).
As \(\|y\|=\|z\|=1\), the Cauchy--Schwarz inequality gives
\[
y^{\Transpose}\PenaltyMatrix{\gamma} z=y^{\Transpose}\OnesMatrix z=\left(\sum_{i\in S}y_{i}\right)\left(\sum_{j\in R}z_{j}\right)\le\sqrt{|S|}\,\sqrt{|R|}\le\sqrt \Biclique .
\]
The first inequality is an equality if and only if \(y\) and \(z\) are uniform on \(S\) and \(R\), and the second if and only if \(|S||R|=\Biclique\).
\end{proof}

Let \(\SmallestPenalty\left(G\right)\coloneqq\inf\left\{\gamma\ge0:\QuadraticMax{\gamma}\left(G\right)=\frac12\sqrt{\Biclique\left(G\right)}\right\}\) be the smallest exact penalty.

\begin{theorem}[]\protect\hypertarget{thm-penalty}{}\label{thm-penalty}

Let \(G\) be bipartite and \(\gamma\ge\frac12\left(1+\sqrt{\Biclique\left(G\right)}\right)\).
Then \(\QuadraticMax{\gamma}\left(G\right)=\tfrac12\sqrt{\Biclique\left(G\right)}\), that is,
\[
\tfrac12\sqrt{\Biclique\left(G\right)}=\max_{u\ge0,\ \|u\|=1}u^{\Transpose}\left(\CostMatrix-\gamma \Adjacency\right)u=\min\left\{\lambda:\lambda \Identity-\CostMatrix+\gamma \Adjacency\in\COP[n]\right\}.
\]
If \(\Biclique\left(G\right)\ge1\), the maximizers are exactly the vectors \(\frac1{\sqrt2}\left(e_{A}+e_{B}\right)\) with \(\left(A,B\right)\) a maximum pair.
In particular, the conclusion holds for every \(\gamma\ge\frac12\left(1+\sqrt{\hat \Biclique}\right)\) with \(\hat \Biclique\ge \Biclique\left(G\right)\), such as \(\gamma\) at least \(\tfrac12+\ConeBicliqueRatio{\DNN}\left(G\right)\), \(\tfrac12+\LasserreBicliqueRatio{1}\left(G\right)\), \(\frac12\left(1+\|\OnesMatrix-\Biadjacency\|_{2}\right)\) or \(\frac14\left(2+n\right)\).

\end{theorem}

\begin{proof}
If \(\Biclique\left(G\right)=0\), then \(G\) is complete bipartite and \(\PenaltyMatrix{\gamma}=\left(1-2\gamma\right)\OnesMatrix\le0\) for \(\gamma\ge\frac12\); hence \(\BilinearMax{\gamma}\le0\) and \(\QuadraticMax{\gamma}=0\) by \hyperref[eq-mu-sigma]{(\ref*{eq-mu-sigma})}.
Assume now \(\Biclique=\Biclique\left(G\right)\ge1\).
By \hyperref[eq-mu-sigma]{(\ref*{eq-mu-sigma})} and \hyperref[eq-mu-cop]{(\ref*{eq-mu-cop})} it suffices to show that \(\BilinearMax{\gamma}=\sqrt \Biclique\).

Whenever \(\BilinearMax{\gamma}>\sqrt \Biclique\), a maximizer \(\left(y,z\right)\) of \hyperref[eq-sigma]{(\ref*{eq-sigma})} must have an edge joining \(\nz\left(y\right)\) and \(\nz\left(z\right)\) by \hyperref[lem-sigma]{Lemma~\ref*{lem-sigma}}, and then \(2\gamma-1<\BilinearMax{\gamma}\le\BilinearMax{0}\) by \hyperref[lem-second-order]{Lemma~\ref*{lem-second-order}}.
Hence \(\BilinearMax{\gamma}=\sqrt \Biclique\) for \(\gamma\) large, and by \hyperref[lem-sigma]{Lemma~\ref*{lem-sigma}} the set \(\left\{\gamma\ge0:\BilinearMax{\gamma}=\sqrt \Biclique\right\}\) is a closed half-line, which by \hyperref[eq-mu-sigma]{(\ref*{eq-mu-sigma})} is \(\left[\SmallestPenalty,\infty\right)\).
It remains to show \(2\SmallestPenalty-1\le\sqrt \Biclique\), which is clear if \(\SmallestPenalty=0\).
Otherwise, for every \(\gamma<\SmallestPenalty\) we have \(\BilinearMax{\gamma}>\sqrt \Biclique\), hence \(2\gamma-1<\BilinearMax{\gamma}\), and letting \(\gamma\) increase to \(\SmallestPenalty\) gives \(2\SmallestPenalty-1\le\BilinearMax{\SmallestPenalty}=\sqrt \Biclique\) by continuity.

Now let \(\gamma\ge\frac12\left(1+\sqrt \Biclique\right)\).
A maximizer of \hyperref[eq-sigma]{(\ref*{eq-sigma})} has value \(\sqrt \Biclique\le2\gamma-1\).
An edge in its support would give \(2\gamma-1<\sqrt \Biclique\) by \hyperref[lem-second-order]{Lemma~\ref*{lem-second-order}}; thus there is none, and by \hyperref[lem-sigma]{Lemma~\ref*{lem-sigma}} it is \(\left(e_{A},e_{B}\right)\) for a maximum pair \(\left(A,B\right)\).
A maximizer of \hyperref[eq-mu]{(\ref*{eq-mu})} has the form \(\left(a\hat y,b\hat z\right)\) with \(ab=\frac12\), that is \(a=b=\frac1{\sqrt2}\), and with \(\left(\hat y,\hat z\right)\) maximizing \hyperref[eq-sigma]{(\ref*{eq-sigma})}, which gives the stated form.
Finally, the four thresholds are at least \(\frac12\left(1+\sqrt \Biclique\right)\).
Indeed, \(\tfrac12\sqrt \Biclique=\ConeBicliqueRatio{\CP}\left(G\right)\le\ConeBicliqueRatio{\DNN}\left(G\right)\le \LasserreBicliqueRatio{1}\left(G\right)\) by \hyperref[thm-cp]{Theorem~\ref*{thm-cp}} and the inclusions \(\CP[n]\subseteq\DNN[n]\subseteq\PSD[n]\); \(\sqrt \Biclique\le\|\OnesMatrix-\Biadjacency\|_{2}\), because \(\OnesMatrix-\Biadjacency\) is the all-ones matrix on \(A\times B\) for a maximum pair \(\left(A,B\right)\); and \(\|\OnesMatrix-\Biadjacency\|_{2}\le\|\OnesMatrix\|_{2}=\sqrt{m_{1}m_{2}}\le\frac n2\), because \(0\le \OnesMatrix-\Biadjacency\le \OnesMatrix\).
The proof of \hyperref[thm-cp]{Theorem~\ref*{thm-cp}} in \hyperref[sec-cp]{Section~\ref*{sec-cp}} uses no result of this section.
\end{proof}

In the parameter \(2\gamma-1\) of \hyperref[eq-sigma]{(\ref*{eq-sigma})}, \hyperref[thm-penalty]{Theorem~\ref*{thm-penalty}} says that a penalty is exact as soon as it is at least the optimal value \(\sqrt{\Biclique}\), whatever the size of \(G\).

\subsection{Computable Thresholds}\label{sec-computable}

The last part of \hyperref[thm-penalty]{Theorem~\ref*{thm-penalty}} replaces the unknown \(\Biclique\left(G\right)\) by an upper bound that can be computed, and the semidefinite bounds listed there are never worse than the spectral one.
If \(X\) is feasible for \hyperref[eq-h1]{(\ref*{eq-h1})}, then \(\left\langle \CostMatrix,X\right\rangle=\left\langle \CostMatrix-\frac12\Adjacency,X\right\rangle\) because \(X_{E}=0\).
The matrix \(\CostMatrix-\frac12\Adjacency\) is the matrix \(\CostMatrix\) built from \(\OnesMatrix-\Biadjacency\) instead of \(\OnesMatrix\), and its largest eigenvalue is \(\frac12\|\OnesMatrix-\Biadjacency\|_{2}\).
Hence
\begin{equation}\protect\phantomsection\label{eq-thresholds}{
\sqrt{\Biclique\left(G\right)}\le2\ConeBicliqueRatio{\DNN}\left(G\right)\le2\LasserreBicliqueRatio{1}\left(G\right)\le\|\OnesMatrix-\Biadjacency\|_{2}\le\|\OnesMatrix-\Biadjacency\|_{F}=\sqrt{\left|E\left(\Complement{G}\right)\right|}.
}\end{equation}
The level-1 semidefinite bound of Parrilo gives a still smaller computable threshold (\hyperref[sec-rho]{Section~\ref*{sec-rho}}).

Because \(\Biclique\) is an integer, an estimate \(\hat \Biclique\ge \Biclique\) in \hyperref[thm-penalty]{Theorem~\ref*{thm-penalty}} can be replaced by \(\lfloor\hat \Biclique\rfloor\), and the thresholds can be rounded, for example to \(\frac12\left(1+\sqrt{\lfloor4\LasserreBicliqueRatio{1}\left(G\right)^{2}\rfloor}\right)\) or \(\frac12\left(1+\sqrt{\lfloor\|\OnesMatrix-\Biadjacency\|_{2}^{2}\rfloor}\right)\).
A semidefinite bound is known only up to the tolerance \(\varepsilon\) of the solver, and we then take \(\gamma\ge\tfrac12+\hat h+\varepsilon\), where \(\hat h\) is the computed value of \(\ConeBicliqueRatio{\DNN}\left(G\right)\) or \(\LasserreBicliqueRatio{1}\left(G\right)\).

\subsection{Sharpness}\label{sec-sharp}

By \hyperref[thm-penalty]{Theorem~\ref*{thm-penalty}} and its proof, the infimum \(\SmallestPenalty\left(G\right)\) is attained and \(2\SmallestPenalty\left(G\right)-1\le\sqrt{\Biclique\left(G\right)}\).
Graphs with a single edge show that the threshold \(\frac12\left(1+\sqrt \Biclique\right)\) of \hyperref[thm-penalty]{Theorem~\ref*{thm-penalty}} cannot be lowered by more than \(\frac12\) in general.

\begin{proposition}[]\protect\hypertarget{prp-sharp}{}\label{prp-sharp}

For an integer \(k\ge1\), let \(N_{k}\) have \(V_{1}=\left\{v_{10},v_{11},\dots,v_{1k}\right\}\), \(V_{2}=\left\{v_{21},\dots,v_{2k}\right\}\) and the single edge \(v_{10}v_{21}\).
Then \(\Biclique\left(N_{k}\right)=k^{2}\) and \(\SmallestPenalty\left(N_{k}\right)=\frac k2=\frac12\sqrt{\Biclique\left(N_{k}\right)}\).

\end{proposition}

\begin{proof}
We write \(y_{i}\) for the entry of \(y\) at \(v_{1i}\) and \(z_{j}\) for the entry of \(z\) at \(v_{2j}\).
A biindependent pair containing \(v_{10}\) avoids \(v_{21}\) and therefore has product at most \(\left(k+1\right)\left(k-1\right)\), while \(\left(\left\{v_{11},\dots,v_{1k}\right\},V_{2}\right)\) has product \(k^{2}\); hence \(\Biclique\left(N_{k}\right)=k^{2}\).
By \hyperref[eq-mu-sigma]{(\ref*{eq-mu-sigma})} we must show that \(\BilinearMax{\gamma}\left(N_{k}\right)>k\) for \(\gamma<\frac k2\) and \(\BilinearMax{\gamma}\left(N_{k}\right)\le k\) for \(\gamma=\frac k2\).

For \(\gamma<\frac k2\), let \(z=e_{V_{2}}\).
Then \(\left(\PenaltyMatrix{\gamma} z\right)_{i}=\sqrt k\) for \(1\le i\le k\) and \(\sqrt k\left(\PenaltyMatrix{\gamma} z\right)_{0}=k-2\gamma>0\); thus \(y=\|\PenaltyMatrix{\gamma} z\|^{-1}\PenaltyMatrix{\gamma} z\ge0\) gives \(\BilinearMax{\gamma}\ge\|\PenaltyMatrix{\gamma} z\|=\sqrt{k^{2}+\frac1k\left(k-2\gamma\right)^{2}}>k\).

For \(\gamma=\frac k2\), the rows of \(\PenaltyMatrix{\gamma}\) indexed by \(v_{11},\dots,v_{1k}\) and the columns indexed by \(v_{22},\dots,v_{2k}\) are all-ones vectors; hence a feasible \(\left(y,z\right)\) enters the objective only through \(y_{0}\), \(z_{1}\) and the sums \(s=\sum_{i\ge1}y_{i}\), \(t=\sum_{j\ge2}z_{j}\):
\[
y^{\Transpose}\PenaltyMatrix{\gamma} z=\left(s+y_{0}\right)\left(t+z_{1}\right)-k\,y_{0}z_{1} .
\]
Write \(y_{0}=\sin\theta\) and \(z_{1}=\sin\omega\) with \(\theta,\omega\in\left[0,\frac\pi2\right]\).
The coefficients of \(s\) and \(t\) are nonnegative, and by the Cauchy--Schwarz inequality \(s\le\sqrt k\cos\theta\) and \(t\le\sqrt{k-1}\cos\omega\).
With \(P=\sqrt k\cos\theta+\sin\theta\ge0\) this gives
\[
y^{\Transpose}\PenaltyMatrix{\gamma} z\le P\sqrt{k-1}\cos\omega+\left(P-k\sin\theta\right)\sin\omega\le\sqrt{\left(k-1\right)P^{2}+\left(P-k\sin\theta\right)^{2}};
\]
the second inequality is Cauchy--Schwarz when \(P-k\sin\theta\ge0\), and otherwise the middle expression is at most \(P\sqrt{k-1}\), which is at most the right side.
Finally, since \(P-\sin\theta=\sqrt k\cos\theta\),
\[
\begin{aligned}
\left(k-1\right)P^{2}+\left(P-k\sin\theta\right)^{2}&=k\left(\left(P-\sin\theta\right)^{2}+\left(k-1\right)\sin^{2}\theta\right)\\
&=k\left(k-\sin^{2}\theta\right)\le k^{2} .
\end{aligned}
\]
Hence \(y^{\Transpose}\PenaltyMatrix{\gamma} z\le k\) for every feasible \(\left(y,z\right)\), that is, \(\BilinearMax{\gamma}\left(N_{k}\right)\le k\).
\end{proof}

The additive \(\frac12\) cannot be dropped altogether, however.
The penalty \(\frac12\sqrt \Biclique\) is exact on \(N_{k}\) by \hyperref[prp-sharp]{Proposition~\ref*{prp-sharp}}, but it is not exact on the graph with \(V_{1}=\left\{v_{11},v_{12}\right\}\), \(V_{2}=\left\{v_{21},v_{22},v_{23}\right\}\) and the edges \(v_{11}v_{22}\) and \(v_{12}v_{21}\), which has \(\Biclique=2\).
For \(y=e_{V_{1}}\) and \(\gamma<1\), the vector \(\PenaltyMatrix{\gamma}^{\Transpose}y=\sqrt2\left(1-\gamma,1-\gamma,1\right)\) is nonnegative with norm \(\sqrt{2+4\left(1-\gamma\right)^{2}}>\sqrt2\).
Hence \(\BilinearMax{\gamma}>\sqrt \Biclique\) for every \(\gamma<1\), and \(\SmallestPenalty\ge1>\frac12\sqrt2=\frac12\sqrt \Biclique\).

\subsection{Local Maximizers and Comparison with Gillis and Glineur}\label{sec-local}

Local search methods, such as the ascent of \hyperref[sec-ascent]{Section~\ref*{sec-ascent}}, reach local maximizers of \hyperref[eq-sigma]{(\ref*{eq-sigma})}, which \hyperref[thm-penalty]{Theorem~\ref*{thm-penalty}} does not describe.

\begin{proposition}[]\protect\hypertarget{prp-local}{}\label{prp-local}

Let \(\Biclique\left(G\right)\ge1\) and \(\gamma\ge\frac12\left(1+\sqrt{\Biclique\left(G\right)}\right)\).
Every local maximizer of \hyperref[eq-sigma]{(\ref*{eq-sigma})} with positive value is \(\left(e_{A},e_{B}\right)\) for a maximal biindependent pair \(\left(A,B\right)\).

\end{proposition}

\begin{proof}
Let \(\left(y,z\right)\) be a local maximizer with value \(\BilinearValue>0\), and let \(S=\nz\left(y\right)\) and \(R=\nz\left(z\right)\).
By \hyperref[thm-penalty]{Theorem~\ref*{thm-penalty}}, \(\BilinearValue\le\BilinearMax{\gamma}=\sqrt \Biclique\le2\gamma-1\); hence by \hyperref[lem-second-order]{Lemma~\ref*{lem-second-order}} no edge joins \(S\) and \(R\), and \(\left(S,R\right)\) is biindependent.
\hyperref[lem-second-order]{Lemma~\ref*{lem-second-order}} gives \(\BilinearValue y_{S}=\left(\PenaltyMatrix{\gamma} z\right)_{S}\), and \(\left(\PenaltyMatrix{\gamma} z\right)_{i}=\OnesVector^{\Transpose}z\) for \(i\in S\) because no edge joins \(i\) and \(R\); thus \(y\) is uniform on \(S\), and likewise \(z\) is uniform on \(R\).
If a vertex \(v\in V_{1}\setminus S\) had no neighbor in \(R\), then \(\left(\PenaltyMatrix{\gamma} z\right)_{v}=\OnesVector^{\Transpose}z>0\), and for \(t\ge0\) the feasible point
\[
y\left(t\right)=\frac{y+te_{v}}{\|y+te_{v}\|}\qquad\text{has value}\qquad y\left(t\right)^{\Transpose}\PenaltyMatrix{\gamma} z=\frac{\BilinearValue+t\,\OnesVector^{\Transpose}z}{\sqrt{1+t^{2}}},
\]
which exceeds \(\BilinearValue\) for small \(t>0\).
This contradicts local maximality, and the same holds for \(V_{2}\setminus R\); hence \(\left(S,R\right)\) is maximal.
\end{proof}

Gillis and Glineur \citep[Section 5]{gillis2014continuous} characterize the maximizers of \hyperref[eq-sigma]{(\ref*{eq-sigma})} under two conditions on the parameter \(2\gamma-1\).
If \(2\gamma-1\ge\sqrt{\left|E\left(\Complement{G}\right)\right|}\), every local maximizer with positive value is a maximal pair, and hence every global maximizer is a maximum pair \citep[Lemma 3]{gillis2014continuous} (see also \citep{gillis2008nonnegative}).
If \(2\gamma-1\ge\max\left(m_{1},m_{2}\right)\), the converse holds as well, that is, every maximal pair is a local maximizer \citep[Theorem 1]{gillis2014continuous}.
Together, \hyperref[thm-penalty]{Theorem~\ref*{thm-penalty}} and \hyperref[prp-local]{Proposition~\ref*{prp-local}} give the first conclusion already for \(2\gamma-1\ge\sqrt{\Biclique\left(G\right)}\).
This threshold is never larger than the two others, because a biindependent pair \(\left(A,B\right)\) spans \(|A||B|\) edges of \(\Complement{G}\), and thus \(\sqrt \Biclique\le\sqrt{\left|E\left(\Complement{G}\right)\right|}\le\sqrt{m_{1}m_{2}}\le\max\left(m_{1},m_{2}\right)\).
By \hyperref[eq-thresholds]{(\ref*{eq-thresholds})}, even the spectral threshold \(2\gamma-1\ge\|\OnesMatrix-\Biadjacency\|_{2}\), which needs one singular value and not \(\Biclique\), is never larger than \(\sqrt{\left|E\left(\Complement{G}\right)\right|}\).
Our threshold does not give the converse, which holds for \(2\gamma-1>\max\left(m_{1},m_{2}\right)-1\) and can fail below this value \citep[Remark 1]{gillis2014continuous}.

\begin{refremark}
At the threshold of \hyperref[thm-penalty]{Theorem~\ref*{thm-penalty}}, the uniform vectors on a maximal pair are not always a local maximizer, and thus the converse of \hyperref[prp-local]{Proposition~\ref*{prp-local}} does not hold.
An example is the graph of Gillis and Glineur \citep[Remark 1]{gillis2014continuous}, with \(V_{1}=\left\{v_{11},v_{12}\right\}\) and \(V_{2}=\left\{v_{21},\dots,v_{24}\right\}\) and the single edge \(v_{12}v_{23}\).
Here \(\Biclique=6\), attained by \(\left(V_{1},V_{2}\setminus\left\{v_{23}\right\}\right)\), and the threshold of \hyperref[thm-penalty]{Theorem~\ref*{thm-penalty}} is \(\frac12\left(1+\sqrt6\right)<2\).
The pair \(\left(\left\{v_{11}\right\},V_{2}\right)\) is maximal, as \(v_{12}\) has a neighbor in \(V_{2}\), but for \(z=e_{V_{2}}\) the entry of \(\PenaltyMatrix{\gamma} z\) at \(v_{12}\) is \(2-\gamma\).
Hence for every \(\gamma<2\), in particular at the threshold, the path \(y\left(t\right)\) of the proof of \hyperref[prp-local]{Proposition~\ref*{prp-local}}, from \(y=e_{v_{11}}\) toward \(v_{12}\), increases the value, and \(\left(e_{v_{11}},e_{V_{2}}\right)\) is not a local maximizer.

\label{rem-converse}

\end{refremark}

\textbf{Pareto singular values.}
Program \hyperref[eq-sigma]{(\ref*{eq-sigma})} is also a Pareto singular value problem in the sense of Seeger and Sossa \citep{seeger2025pareto}.
The value of \hyperref[eq-sigma]{(\ref*{eq-sigma})} is, up to sign, the least Pareto singular value of \(-\left(\OnesMatrix-2\gamma \Biadjacency\right)\), a quantity that is NP-hard to compute for general matrices \citep[Theorem 2]{barbarino2026computing}.
For a nonnegative matrix the maximum in \hyperref[eq-sigma]{(\ref*{eq-sigma})} is the spectral norm, by the Perron--Frobenius theorem, whereas here the negative entries \(1-2\gamma\) on the edges of \(G\) encode biindependence.
Seeger and Sossa study the Pareto singular values of Boolean matrices and their relation to the spectra of bipartite graphs, but they do not consider penalized matrices or penalty thresholds.

\section{Exact Copositive and Completely Positive Formulations}\label{sec-cp}

Copositive and completely positive programs give exact convex formulations of many nonconvex quadratic and combinatorial problems \citep{burer2009copositive, bomze2012copositive}.
For the stability number \(\Independence\left(G\right)\), the largest size of a stable set of \(G\), replacing the positive semidefinite cone in the Lovász theta number \(\LovaszTheta\left(G\right)\) \citep{lovasz1979shannon} by the completely positive cone gives \(\Independence\left(G\right)\) itself \citep[Theorem 2.2]{deklerk2002approximation}.
The dual statement is \(\Independence\left(G\right)=\min\left\{\lambda:\lambda\left(\Identity+\Adjacency\right)-\OnesMatrix[n]\ \text{copositive}\right\}\), where \(\Adjacency\) is the adjacency matrix of \(G\) and \(\OnesMatrix[n]\) the \(n\times n\) all-ones matrix, and it is a version of the Motzkin--Straus theorem \citep[Corollary 2.4]{deklerk2002approximation}.
We show that \(\tfrac12\sqrt{\Biclique}\) has the same two descriptions, with the bound \(\LasserreBicliqueRatio{1}\) of \hyperref[eq-h1]{(\ref*{eq-h1})} in place of the theta number.
The completely positive version of \(\LasserreBicliqueRatio{1}\) is exact (\hyperref[thm-cp]{Theorem~\ref*{thm-cp}}), and its dual contains the penalized program of \hyperref[thm-penalty]{Theorem~\ref*{thm-penalty}}, which makes both optima attained (\hyperref[cor-dual]{Corollary~\ref*{cor-dual}}).

\begin{theorem}[]\protect\hypertarget{thm-cp}{}\label{thm-cp}

For every bipartite graph \(G\), \(\ConeBicliqueRatio{\CP}\left(G\right)=\tfrac12\sqrt{\Biclique\left(G\right)}\), and the maximum is attained.

\end{theorem}

\begin{proof}
The feasible set is nonempty (it contains \(e_{i}e_{i}^{\Transpose}\)), closed (\(\CP[n]\) is closed) and bounded (\(\tr\left(X\right)=1\)); hence the maximum is attained.

Let \(X\) be feasible and write \(X=\sum_{t}\lambda_{t}u_{t}u_{t}^{\Transpose}\) with \(\lambda_{t}>0\), \(u_{t}\ge0\) and \(\|u_{t}\|=1\); then \(\sum_{t}\lambda_{t}=\tr\left(X\right)=1\).
For \(ij\in E\), \(0=X_{ij}=\sum_{t}\lambda_{t}\left(u_{t}\right)_{i}\left(u_{t}\right)_{j}\) is a sum of nonnegative terms; thus \(\left(u_{t}\right)_{i}\left(u_{t}\right)_{j}=0\) for every \(t\).
Hence \(\nz\left(u_{t}\right)\) is stable, and \(A_{t}=\nz\left(u_{t}\right)\cap V_{1}\), \(B_{t}=\nz\left(u_{t}\right)\cap V_{2}\) form a biindependent pair.
With \(a_{t}=\|u_{t}|_{A_{t}}\|\) and \(b_{t}=\|u_{t}|_{B_{t}}\|\) we have \(a_{t}^{2}+b_{t}^{2}=1\), hence \(a_{t}b_{t}\le\frac12\), and the Cauchy--Schwarz inequality gives
\[
\begin{aligned}
u_{t}^{\Transpose}\CostMatrix u_{t}&=\left(\sum_{i\in A_{t}}\left(u_{t}\right)_{i}\right)\left(\sum_{j\in B_{t}}\left(u_{t}\right)_{j}\right)\\
&\le\sqrt{|A_{t}|}\,a_{t}\sqrt{|B_{t}|}\,b_{t}\le\tfrac12\sqrt{|A_{t}||B_{t}|}\le\tfrac12\sqrt{\Biclique\left(G\right)} .
\end{aligned}
\]
This also holds when \(A_{t}\) or \(B_{t}\) is empty, as the left side is then \(0\).
Summing over \(t\) gives \(\left\langle \CostMatrix,X\right\rangle\le\tfrac12\sqrt{\Biclique\left(G\right)}\).

Conversely, if \(\Biclique\left(G\right)=0\), the value \(0\) is attained by \(e_{i}e_{i}^{\Transpose}\).
Otherwise let \(\left(A,B\right)\) be a maximum pair and \(u=\frac1{\sqrt2}\left(e_{A}+e_{B}\right)\).
Then \(\|u\|^{2}=1\), \(X=uu^{\Transpose}\) lies in \(\CP[n]\) and vanishes on \(E\) because no edge joins \(A\) and \(B\), and \(\left\langle \CostMatrix,X\right\rangle=u^{\Transpose}\CostMatrix u\) is the product of the sums of \(u\) over \(A\) and over \(B\), which is \(\tfrac12\sqrt{|A||B|}\).
\end{proof}

The proof follows the extreme-ray argument of de Klerk and Pasechnik \citep[Theorem 2.2]{deklerk2002approximation} for \(\Independence\left(G\right)\).
Since \(\left\langle \CostMatrix,X\right\rangle\) couples only the two sides, an extreme ray \(uu^{\Transpose}\) is here bounded by the Cauchy--Schwarz inequality on each side.

\begin{refremark}
If \(\Biclique\left(G\right)\ge1\), equality in the proof of \hyperref[thm-cp]{Theorem~\ref*{thm-cp}} forces every \(u_{t}\) to be of the form \(\frac1{\sqrt2}\left(e_{A_{t}}+e_{B_{t}}\right)\) with \(|A_{t}||B_{t}|=\Biclique\left(G\right)\).
The optimal solutions are therefore exactly the convex combinations of the matrices \(uu^{\Transpose}\) for such vectors.

\label{rem-optimal}

\end{refremark}

The copositive program dual to \(\ConeBicliqueRatio{\CP}\) has one multiplier for each edge, and the penalty of \hyperref[thm-penalty]{Theorem~\ref*{thm-penalty}} gives an optimal solution of it.

\begin{corollary}[]\protect\hypertarget{cor-dual}{}\label{cor-dual}

Let \(S_{ij}\coloneqq e_{i}e_{j}^{\Transpose}+e_{j}e_{i}^{\Transpose}\) for \(ij\in E\).
The program \(\ConeBicliqueRatio{\CP}\left(G\right)\) and
\begin{equation}\protect\phantomsection\label{eq-dual}{
\min\left\{\lambda:\lambda \Identity-\CostMatrix+\sum_{ij\in E}\EdgeMultiplier_{ij}S_{ij}\in\COP[n],\ \lambda\in\Reals,\ \EdgeMultiplier\in\Reals^{E}\right\}
}\end{equation}
have the same value \(\tfrac12\sqrt{\Biclique\left(G\right)}\), and both are attained.
An optimal solution of \hyperref[eq-dual]{(\ref*{eq-dual})} is \(\lambda=\tfrac12\sqrt{\Biclique\left(G\right)}\) and \(\EdgeMultiplier_{ij}=\gamma\) for all \(ij\in E\), for any \(\gamma\ge\frac12\left(1+\sqrt{\Biclique\left(G\right)}\right)\).

\end{corollary}

\begin{proof}
For feasible \(X\) and \(\left(\lambda,\EdgeMultiplier\right)\), and with \(Y=\sum_{ij\in E}\EdgeMultiplier_{ij}S_{ij}\), we have \(\left\langle\lambda \Identity-\CostMatrix+Y,X\right\rangle\ge0\) because \(\CP[n]\) and \(\COP[n]\) are dual, and since \(\left\langle Y,X\right\rangle=0\) and \(\tr\left(X\right)=1\) this gives \(\lambda\ge\left\langle \CostMatrix,X\right\rangle\).
The primal value \(\tfrac12\sqrt \Biclique\) is attained by \hyperref[thm-cp]{Theorem~\ref*{thm-cp}}, and the stated point is feasible for \hyperref[eq-dual]{(\ref*{eq-dual})} with the same value by \hyperref[thm-penalty]{Theorem~\ref*{thm-penalty}}.
\end{proof}

The role of \hyperref[thm-penalty]{Theorem~\ref*{thm-penalty}} in \hyperref[cor-dual]{Corollary~\ref*{cor-dual}} is to supply the attainment of the dual optimum.
Program \hyperref[eq-dual]{(\ref*{eq-dual})} is strictly feasible, since \(\lambda \Identity-\CostMatrix\) is positive definite, hence strictly copositive, for \(\lambda>\frac12\sqrt{m_{1}m_{2}}\), and conic duality therefore gives equal values and primal attainment.
It does not give attainment in \hyperref[eq-dual]{(\ref*{eq-dual})}, because the primal is never strictly feasible when \(E\neq\varnothing\), as interior points of \(\CP[n]\) have positive entries \citep[Section 3]{dur2010copositive}.

In the copositive formulation of \(\Independence\left(G\right)\) the edges carry the same multiplier \(\lambda\) as the diagonal, and exactness is the Motzkin--Straus theorem.
For \(\tfrac12\sqrt \Biclique\) the edge multipliers are separate variables.
The statement of \hyperref[cor-dual]{Corollary~\ref*{cor-dual}} that a uniform multiplier of at least \(\frac12\left(1+\sqrt \Biclique\right)\) is optimal has no counterpart in the work of de Klerk and Pasechnik \citep{deklerk2002approximation} or of Peña et al. \citep{pena2007computing}.

\section{Tractable Relaxations}\label{sec-relax}

The exact programs of \hyperref[sec-cp]{Section~\ref*{sec-cp}} are intractable in general, whereas \hyperref[thm-penalty]{Theorem~\ref*{thm-penalty}} needs only an upper bound on \(\Biclique\), and this section asks which tractable bound to compute for it.
The standard approach replaces the completely positive cone by a tractable outer approximation, or the copositive cone by tractable inner approximations \citep{deklerk2002approximation, parrilo2000structured}, which for \(\Independence\left(G\right)\) give hierarchies of linear and semidefinite bounds \citep{deklerk2002approximation, parrilo2000structured, pena2007computing}.
Each such relaxation gives a bound of at least \(\tfrac12\sqrt{\Biclique}\), and hence an exact penalty by \hyperref[thm-penalty]{Theorem~\ref*{thm-penalty}}, which is smaller when the bound is closer to \(\tfrac12\sqrt{\Biclique}\).
Our answer is the level-1 semidefinite bound of Parrilo, which is never weaker than the doubly nonnegative bound and is exact on \(\Value{rhooneexact}\) of the \(\Value{graphs}\) random graphs of \hyperref[sec-bounds]{Section~\ref*{sec-bounds}} (\hyperref[sec-rho]{Section~\ref*{sec-rho}}).

\subsection{Doubly Nonnegative Relaxation}\label{sec-dnn}

The first relaxation replaces the completely positive cone in \hyperref[thm-cp]{Theorem~\ref*{thm-cp}} by the doubly nonnegative cone, which gives the bound \(\ConeBicliqueRatio{\DNN}\).
Up to a factor \(\tfrac12\) it is the bound studied by Suda and Tanaka \citep{suda2014cross} and Suda et al. \citep{suda2017semidefinite}, which in the form recorded by Laurent et al. \citep[Remark 2]{laurent2025semidefinite} is
\begin{equation}\protect\phantomsection\label{eq-T}{
\SudaTanakaBound\left(G\right)\coloneqq\max\left\{\left\langle \CostMatrix,X\right\rangle:X\in\DNN[n],\ \tr\left(X_{V_{1}}\right)=\tr\left(X_{V_{2}}\right)=1,\ X_{E}=0\right\},
}\end{equation}
where \(X_{V_{i}}\) is the principal submatrix of \(X\) on \(V_{i}\).
It satisfies \(\sqrt \Biclique\le \SudaTanakaBound\le2\LasserreBicliqueRatio{1}\), and without the constraint \(X\ge0\) its value is \(2\LasserreBicliqueRatio{1}\) \citep[Remark 2]{laurent2025semidefinite}.
The two normalizations of the trace give the same bound, up to the factor \(\frac12\), for every cone that is invariant under scaling the two sides.

\begin{proposition}[]\protect\hypertarget{prp-scaling}{}\label{prp-scaling}

Let \(\GenericCone\subseteq\PSD[n]\) be a closed convex cone that contains \(e_{i}e_{i}^{\Transpose}+e_{j}e_{j}^{\Transpose}\) for some \(i\in V_{1}\), \(j\in V_{2}\), and with \(\Delta X\Delta\in\GenericCone\) for every \(X\in\GenericCone\) and every \(\Delta=\Diag\left(a\OnesVector[V_{1}]+b\OnesVector[V_{2}]\right)\), \(a,b>0\).
Then
\[
\ConeBicliqueRatio{\GenericCone}\left(G\right)=\tfrac12\max\left\{\left\langle \CostMatrix,X\right\rangle:X\in\GenericCone,\ \tr\left(X_{V_{1}}\right)=\tr\left(X_{V_{2}}\right)=1,\ X_{E}=0\right\}.
\]
In particular \(\ConeBicliqueRatio{\DNN}\left(G\right)=\frac12\SudaTanakaBound\left(G\right)\), and \(\frac12\SudaTanakaBound^{\CP}\left(G\right)=\ConeBicliqueRatio{\CP}\left(G\right)=\frac12\sqrt{\Biclique\left(G\right)}\), where \(\SudaTanakaBound^{\CP}\left(G\right)\) is defined as in \hyperref[eq-T]{(\ref*{eq-T})} with \(\CP[n]\) in place of \(\DNN[n]\).

\end{proposition}

\begin{proof}
Let \(\tau\) be the maximum on the right; its feasible set is compact and contains \(e_{i}e_{i}^{\Transpose}+e_{j}e_{j}^{\Transpose}\), of value \(0\); thus \(\tau\ge0\) is attained.
If \(X\) is feasible for \(\tau\), then \(\frac12X\) is feasible for \hyperref[eq-h1K]{(\ref*{eq-h1K})}; hence \(\ConeBicliqueRatio{\GenericCone}\left(G\right)\ge\frac12\tau\).
Conversely, let \(X\) be feasible for \hyperref[eq-h1K]{(\ref*{eq-h1K})} with \(\left\langle \CostMatrix,X\right\rangle>0\), and let \(t_{i}=\tr\left(X_{V_{i}}\right)\), so that \(t_{1}+t_{2}=1\).
Both \(t_{i}\) are positive, since \(t_{1}=0\) would make the \(V_{1}\) rows of \(X\succeq0\) vanish, and then \(\left\langle \CostMatrix,X\right\rangle=0\).
With \(a=t_{1}^{-1/2}\) and \(b=t_{2}^{-1/2}\), the matrix \(\Delta X\Delta\) is feasible for \(\tau\) and has value \(ab\left\langle \CostMatrix,X\right\rangle\); thus \(\tau\ge ab\left\langle \CostMatrix,X\right\rangle\ge2\left\langle \CostMatrix,X\right\rangle\), as \(ab=\left(t_{1}t_{2}\right)^{-1/2}\) and \(t_{1}t_{2}\le\frac14\).
The cones \(\DNN[n]\) and \(\CP[n]\) satisfy the hypotheses, and \hyperref[thm-cp]{Theorem~\ref*{thm-cp}} gives the last equality.
\end{proof}

By \hyperref[prp-scaling]{Proposition~\ref*{prp-scaling}}, the inequalities \(\sqrt \Biclique\le \SudaTanakaBound\le2\LasserreBicliqueRatio{1}\) read \(\ConeBicliqueRatio{\CP}\left(G\right)\le \ConeBicliqueRatio{\DNN}\left(G\right)\le \LasserreBicliqueRatio{1}\left(G\right)\).
They are the analogue of \(\Independence\le\SchrijverTheta\le\LovaszTheta\), where the bound \(\SchrijverTheta\) of Schrijver \citep{schrijver1979comparison} adds the constraint \(X\ge0\) to \(\LovaszTheta\).
The first inequality is an equality on small graphs.

\begin{proposition}[]\protect\hypertarget{prp-small}{}\label{prp-small}

If \(n\le5\), then \(\ConeBicliqueRatio{\DNN}\left(G\right)=\frac12\sqrt{\Biclique\left(G\right)}\).

\end{proposition}

\begin{proof}
For \(n\le4\) this follows from \hyperref[thm-cp]{Theorem~\ref*{thm-cp}}, since \(\DNN[n]=\CP[n]\) for \(n\le4\) \citep{maxfield1962matrix}.
The case \(n=5\) is proved in \ref{sec-five}.
\end{proof}

\begin{refremark}
The result of Kim et al. \citep{kim2020doubly} gives equality of doubly nonnegative and completely positive relaxations when the aggregate sparsity pattern is a block-clique graph with cliques of at most four nodes, and one may ask whether it implies \hyperref[prp-small]{Proposition~\ref*{prp-small}}.
It does not apply here when \(m_{1},m_{2}\ge2\).
The aggregate sparsity pattern of \hyperref[eq-h1K]{(\ref*{eq-h1K})} is the support of \(\CostMatrix\), which is the complete bipartite graph \(K_{m_{1},m_{2}}\), since the constraints \(X_{E}=0\) only involve pairs in this support and the trace constraint only the diagonal.
When \(m_{1},m_{2}\ge2\), this graph contains a \(4\)-cycle without a chord.
It is therefore not chordal, and hence not a block-clique graph.

\label{rem-small}

\end{refremark}

The first inequality can be strict from six vertices on.

\begin{proposition}[]\protect\hypertarget{prp-matching}{}\label{prp-matching}

For \(s\ge2\), let \(sK_{2}\) be the perfect matching on \(s+s\) vertices, made of \(s\) disjoint edges.
Then \(\ConeBicliqueRatio{\DNN}\left(sK_{2}\right)=\LasserreBicliqueRatio{1}\left(sK_{2}\right)=\frac s4\), while \(\tfrac12\sqrt{\Biclique\left(sK_{2}\right)}=\tfrac14\sqrt{s^{2}-1}\) for \(s\) odd.
In particular \(\ConeBicliqueRatio{\DNN}\left(3K_{2}\right)=\frac34>\frac12\sqrt2=\ConeBicliqueRatio{\CP}\left(3K_{2}\right)\).

\end{proposition}

\begin{proof}
Order \(V_{2}\) so that the biadjacency matrix is \(\Identity[s]\).
A pair \(\left(A,B\right)\) is biindependent if and only if \(B\) avoids the partners of \(A\); thus \(\Biclique\left(sK_{2}\right)=\max_{a} a\left(s-a\right)=\lfloor\frac s2\rfloor\lceil\frac s2\rceil\).
For the upper bound \(\LasserreBicliqueRatio{1}\left(sK_{2}\right)\le \frac s4\) \citep[Section 5]{laurent2025semidefinite}, let \(Y=\frac s4\Adjacency\).
The matrix \(\frac s4\Identity-\CostMatrix+Y\) has diagonal blocks \(\frac s4\Identity[s]\) and off-diagonal blocks \(L=\frac s4\Identity[s]-\frac12\OnesMatrix[s]\), whose eigenvalues are \(\frac s4\) on \(\OnesVector^{\perp}\) and \(-\frac s4\) on \(\OnesVector\); as the eigenvalues of the block matrix with diagonal blocks \(aI\) and off-diagonal blocks \(L\) are \(a\pm\lambda\) for the eigenvalues \(\lambda\) of \(L\), it is positive semidefinite.
If \(X\) is feasible for \hyperref[eq-h1]{(\ref*{eq-h1})}, then \(\left\langle Y,X\right\rangle=0\); hence \(\left\langle \CostMatrix,X\right\rangle=\frac s4-\left\langle\frac s4\Identity-\CostMatrix+Y,X\right\rangle\le\frac s4\).
Since \(\ConeBicliqueRatio{\DNN}\left(sK_{2}\right)\le \LasserreBicliqueRatio{1}\left(sK_{2}\right)\), it remains to exhibit a feasible \(X\in\DNN[2s]\) with value \(\frac s4\).
Let
\[
X=\begin{pmatrix}P&Q\\Q&P\end{pmatrix},\qquad P=p\Identity[s]+q\left(\OnesMatrix[s]-\Identity[s]\right),\qquad Q=t\left(\OnesMatrix[s]-\Identity[s]\right),
\]
with
\[
p=\frac1{2s},\qquad t=\frac1{4\left(s-1\right)},\qquad q=p-t=\frac{s-2}{4s\left(s-1\right)}.
\]
Then \(X\ge0\), \(\tr\left(X\right)=2sp=1\), \(X\) vanishes on the matching edges (the diagonal of \(Q\)), and \(\left\langle \CostMatrix,X\right\rangle=\OnesVector^{\Transpose}Q\OnesVector=st\left(s-1\right)=\frac s4\).
The eigenvalues of \(X\) are those of \(P\pm Q=\left(p-q\mp t\right)\Identity[s]+\left(q\pm t\right)\OnesMatrix[s]\), namely \(p+\left(s-1\right)\left(q\pm t\right)\) on \(\OnesVector\) and \(p-q\mp t\) on \(\OnesVector^{\perp}\).
With \(q=p-t\) these are \(sp=\frac12\) and \(0\) for \(P+Q\), and \(sp-2t\left(s-1\right)=0\) and \(2t\) for \(P-Q\); thus \(X\succeq0\).
\end{proof}

By \hyperref[prp-small]{Proposition~\ref*{prp-small}}, \(3K_{2}\) is a smallest graph on which \(\ConeBicliqueRatio{\DNN}\left(G\right)>\frac12\sqrt \Biclique\).

Laurent et al. \citep{laurent2025semidefinite} introduced the bound \(\LasserreBicliqueRatio{1}\) for the ratio \(\BicliqueRatio\) of \hyperref[sec-prelim]{Section~\ref*{sec-prelim}}, which can be smaller than \(\tfrac12\sqrt{\Biclique}\).

\begin{refremark}
Since every cone \(\GenericCone\supseteq\CP[n]\) gives \(\ConeBicliqueRatio{\GenericCone}\ge\tfrac12\sqrt \Biclique\), no bound of the form \hyperref[eq-h1K]{(\ref*{eq-h1K})} can reach \(\BicliqueRatio\) when \(\BicliqueRatio<\tfrac12\sqrt \Biclique\).
This happens exactly when no maximum pair is balanced.
To see it, let \(\BalancedBiclique\left(G\right)\) be the largest product \(|A||B|\) over the balanced biindependent pairs, those with \(|A|=|B|\), as in Laurent et al. \citep{laurent2025semidefinite}.
Then \(\BicliqueRatio\left(G\right)=\tfrac12\sqrt{\Biclique\left(G\right)}\) if and only if \(\BalancedBiclique\left(G\right)=\Biclique\left(G\right)\).
In one direction, a balanced maximum pair gives \(\BicliqueRatio\ge\frac12|A|=\tfrac12\sqrt \Biclique\).
In the other, if \(\left(A,B\right)\) attains \(\BicliqueRatio\left(G\right)\), then \(\BicliqueRatio\left(|A|+|B|\right)=|A||B|\), and \(|A|+|B|\ge2\sqrt{|A||B|}\) gives \(\BicliqueRatio\le\tfrac12\sqrt{|A||B|}\le\tfrac12\sqrt \Biclique\) \citep[Lemma 1]{laurent2025semidefinite}, with equality only if \(|A|=|B|\) and \(|A||B|=\Biclique\).
Two examples have \(\BicliqueRatio<\tfrac12\sqrt \Biclique\).
The graph of Laurent et al. \citep[Example 1]{laurent2025semidefinite}, with one edge on \(2+2\) vertices, has \(\Biclique=2\) and \(\BalancedBiclique=1\), and hence \(\BicliqueRatio=\tfrac23<\tfrac12\sqrt2=\tfrac12\sqrt \Biclique\).
The perfect matchings \(sK_{2}\) of \hyperref[prp-matching]{Proposition~\ref*{prp-matching}} with odd \(s\) also have \(\BicliqueRatio<\tfrac12\sqrt \Biclique\), because \(\Biclique=\frac14\left(s^{2}-1\right)\) is attained only by pairs with sides of sizes \(\frac12\left(s-1\right)\) and \(\frac12\left(s+1\right)\).
Two other examples have \(\BalancedBiclique=\Biclique\) and \(\BicliqueRatio=\tfrac12\sqrt \Biclique\).
On \(sK_{2}\) with even \(s\), a pair with both sides of size \(\frac s2\) attains \(\Biclique\).
On the crown \(K_{s,s}\) minus a perfect matching, the biindependent pairs with both sides nonempty are the matched pairs of single vertices.

\label{rem-ratio}

\end{refremark}

\subsection{The Level-1 Semidefinite Bound}\label{sec-rho}

The second relaxation acts on the dual side, where Parrilo \citep{parrilo2000structured} approximates the copositive cone from inside by cones defined through sums of squares.
The level-1 cone, \(\PolyaSOSCone[n]{1}\), is described by the linear matrix inequalities (13)--(16) of de Klerk and Pasechnik \citep{deklerk2002approximation}, which we recall in \ref{sec-comp}.
Replacing \(\COP[n]\) by \(\PolyaSOSCone[n]{1}\) in \hyperref[eq-dual]{(\ref*{eq-dual})} gives the semidefinite bound
\[
\PolyaSOS{1}\left(G\right)\coloneqq\min\left\{\lambda:\lambda \Identity-\CostMatrix+\sum_{ij\in E}\EdgeMultiplier_{ij}S_{ij}\in\PolyaSOSCone[n]{1}\right\}.
\]
The same program over \(\PSD[n]+\NN[n]\), the dual cone of \(\DNN[n]\), is strictly feasible, and thus it has the same value as \(\ConeBicliqueRatio{\DNN}\) \citep[Theorem 2.4.1]{bental2001lectures}.
As \(\PSD[n]+\NN[n]\subseteq\PolyaSOSCone[n]{1}\subseteq\COP[n]\), it follows that \(\tfrac12\sqrt \Biclique\le\PolyaSOS{1}\le\ConeBicliqueRatio{\DNN}\le \LasserreBicliqueRatio{1}\).
Hence \(\gamma\ge\tfrac12+\PolyaSOS{1}\left(G\right)\) is a computable threshold for \hyperref[thm-penalty]{Theorem~\ref*{thm-penalty}}, never larger than those based on \(\ConeBicliqueRatio{\DNN}\) and \(\LasserreBicliqueRatio{1}\), and it equals \(\frac12\left(1+\sqrt \Biclique\right)\) whenever \(\PolyaSOS{1}\) is exact.

\subsection{Pólya Approximations}\label{sec-polya}

The third relaxation uses the simplest inner approximations of \(\COP[n]\), which are polyhedral and give linear programming bounds.
The penalty they give always exceeds the threshold of \hyperref[thm-penalty]{Theorem~\ref*{thm-penalty}} (\hyperref[thm-polya-inexact]{Theorem~\ref*{thm-polya-inexact}}).
For \(M\in\Sym[n]\) we write \(p_{M}\left(\FormVariable\right)\coloneqq\sum_{i,j}M_{ij}\FormVariable_{i}^{2}\FormVariable_{j}^{2}\), so that \(M\in\COP[n]\) if and only if \(p_{M}\ge0\) on \(\Reals^{n}\), as \(p_{M}\left(\FormVariable\right)=u^{\Transpose}Mu\) for \(u=\left(\FormVariable_{1}^{2},\dots,\FormVariable_{n}^{2}\right)\).
For \(\nu\in\Naturals^{n}\), with \(\Naturals\coloneqq\left\{0,1,2,\dots\right\}\), we write \(\FormVariable^{\nu}\coloneqq\prod_{i}\FormVariable_{i}^{\nu_{i}}\) and \(\nu!\coloneqq\prod_{i}\nu_{i}!\), and \(\Naturals^{n}_{d}\coloneqq\left\{\nu\in\Naturals^{n}:\sum_{i}\nu_{i}=d\right\}\).
With \(\NonnegativeForms\) the set of polynomials in \(\FormVariable\in\Reals^{n}\) with nonnegative coefficients, the Pólya cones of de Klerk and Pasechnik \citep[Definition 3.1]{deklerk2002approximation} are
\[
\PolyaCone[n]{r}\coloneqq\left\{M\in\Sym[n]:\left(\sum_{i}\FormVariable_{i}^{2}\right)^{r}p_{M}\in\NonnegativeForms\right\},\qquad r\ge0,
\]
and they satisfy \(\NN[n]=\PolyaCone[n]{0}\subseteq\PolyaCone[n]{1}\subseteq\dots\subseteq\COP[n]\).
Expanded, \(\left(\sum_{i}\FormVariable_{i}^{2}\right)^{r}p_{M}\) has only the monomials \(\FormVariable^{2\nu}\), for \(\nu\in\Naturals^{n}_{r+2}\), and the coefficient of \(\FormVariable^{2\nu}\) is \(\frac{r!}{\nu!}\left(\nu^{\Transpose}M\nu-\nu^{\Transpose}\diag M\right)\).
Thus \(M\in\PolyaCone[n]{r}\) if and only if \(\nu^{\Transpose}M\nu-\nu^{\Transpose}\diag M\ge0\) for all such \(\nu\) \citep{bomze2002solving}, which is a set of linear inequalities in \(M\).
We put
\begin{equation}\protect\phantomsection\label{eq-zeta}{
\Polya{r}\left(G\right)\coloneqq\min\left\{\lambda:\lambda \Identity-\CostMatrix+\sum_{ij\in E}\EdgeMultiplier_{ij}S_{ij}\in\PolyaCone[n]{r},\ \lambda\in\Reals,\ \EdgeMultiplier\in\Reals^{E}\right\},
}\end{equation}
with \(\Polya{r}\left(G\right)=\infty\) if the linear program is infeasible.
Then \(\Polya{r}\left(G\right)\ge\tfrac12\sqrt{\Biclique\left(G\right)}\), and \(\Polya{r}\left(G\right)\) is nonincreasing in \(r\).
These bounds converge, because the matrix \(\lambda \Identity-\CostMatrix+\gamma \Adjacency\) of \hyperref[thm-penalty]{Theorem~\ref*{thm-penalty}} is strictly copositive for every \(\lambda>\tfrac12\sqrt{\Biclique}\), and Pólya's theorem applies to it.

\begin{proposition}[]\protect\hypertarget{prp-polya-converge}{}\label{prp-polya-converge}

If \(\Biclique\left(G\right)\ge1\), then \(\Polya{r}\left(G\right)\) decreases to \(\tfrac12\sqrt{\Biclique\left(G\right)}\) as \(r\to\infty\).

\end{proposition}

\begin{proof}
Fix any \(\lambda>\tfrac12\sqrt \Biclique\), and let \(\gamma\ge\frac12\left(1+\sqrt \Biclique\right)\).
By \hyperref[thm-penalty]{Theorem~\ref*{thm-penalty}}, \(u^{\Transpose}\left(\lambda \Identity-\CostMatrix+\gamma \Adjacency\right)u\ge\lambda-\tfrac12\sqrt \Biclique>0\) for every unit vector \(u\ge0\); hence \(\lambda \Identity-\CostMatrix+\gamma \Adjacency\) is strictly copositive.
By Pólya's theorem \citep{polya1928positive}, in the form \citep[Theorem 3.3]{deklerk2002approximation}, it lies in \(\PolyaCone[n]{r}\) for all large \(r\); thus \(\Polya{r}\left(G\right)\le\lambda\) for those \(r\).
\end{proof}

Convergence says nothing about the first level at which \hyperref[eq-zeta]{(\ref*{eq-zeta})} is feasible.
This level depends on
\[
\Biindependence\left(G\right)\coloneqq\max\left\{|A|+|B|:\left(A,B\right)\ \text{biindependent},\ A\ne\varnothing,\ B\ne\varnothing\right\},
\]
which is defined when \(\Biclique\left(G\right)\ge1\) and satisfies \(2\le\Biindependence\left(G\right)\le\Independence\left(G\right)\).

\begin{theorem}[]\protect\hypertarget{thm-polya-finite}{}\label{thm-polya-finite}

Let \(\Biclique\left(G\right)\ge1\).
Then \(\Polya{r}\left(G\right)<\infty\) if and only if \(r\ge\Biindependence\left(G\right)-1\).

\end{theorem}

\begin{proof}
Let \(M=\lambda \Identity-\CostMatrix+Y\) with \(Y=\sum_{ij\in E}\EdgeMultiplier_{ij}S_{ij}\).
By the criterion of \citep{bomze2002solving} recalled above, \(M\in\PolyaCone[n]{r}\) if and only if \(q\left(\nu\right)=\nu^{\Transpose}M\nu-\nu^{\Transpose}\diag M\ge0\) for every \(\nu\in\Naturals^{n}_{r+2}\).
As the diagonal of \(M\) is \(\lambda\OnesVector\),
\begin{equation}\protect\phantomsection\label{eq-polya-coeff}{
q\left(\nu\right)=\underbrace{\lambda\sum_{i}\nu_{i}\left(\nu_{i}-1\right)}_{\text{diagonal term}}-\underbrace{\left(\sum_{i\in V_{1}}\nu_{i}\right)\left(\sum_{j\in V_{2}}\nu_{j}\right)}_{\text{cross term}}+\underbrace{2\sum_{ij\in E}\EdgeMultiplier_{ij}\nu_{i}\nu_{j}}_{\text{edge term}} .
}\end{equation}
For \(\lambda\ge0\), the diagonal term is nonnegative and vanishes when \(\nu\) is a \(0/1\) vector; the cross term vanishes when \(\nz\left(\nu\right)\) lies in one side; and the edge term vanishes when no edge joins two vertices of \(\nz\left(\nu\right)\).

Suppose first that \(r+2\le\Biindependence\), and let \(\left(A,B\right)\) be a biindependent pair with both sides nonempty and \(|A|+|B|=\Biindependence\).
Take a \(0/1\) vector \(\nu\) with \(r+2\) entries equal to \(1\), all in \(A\cup B\), and at least one in each of \(A\) and \(B\).
Its diagonal and edge terms vanish, and its cross term is at least \(1\); hence \(q\left(\nu\right)<0\) whatever \(\lambda\) and \(\EdgeMultiplier\) are, and the program \hyperref[eq-zeta]{(\ref*{eq-zeta})} is infeasible.
This is the argument of \citep[Theorem 4.2]{deklerk2002approximation} for \(\Independence\left(G\right)\).

Suppose now that \(r+2>\Biindependence\), and take \(\lambda=t\) and \(\EdgeMultiplier_{ij}=t\) for every edge, with \(t\ge\frac18\left(r+2\right)^{2}\).
Let \(\nu\in\Naturals^{n}_{r+2}\).
The diagonal and edge terms of \hyperref[eq-polya-coeff]{(\ref*{eq-polya-coeff})} are then nonnegative, and the cross term is at most \(\frac14\left(r+2\right)^{2}\le2t\); thus it suffices that one of the other two terms is at least \(2t\) or that the cross term vanishes.
If some edge \(ij\) joins two vertices of \(\nz\left(\nu\right)\), the edge term is at least \(2t\nu_{i}\nu_{j}\ge2t\).
If \(\nz\left(\nu\right)\) lies in one side, the cross term vanishes.
Otherwise \(\nz\left(\nu\right)\) meets both sides and contains no edge, and thus it is a biindependent pair with both sides nonempty, and \(|{\nz\left(\nu\right)}|\le\Biindependence<r+2=\sum_{i}\nu_{i}\); hence some \(\nu_{i}\ge2\), and the diagonal term is at least \(2t\).
In each case \(q\left(\nu\right)\ge0\); thus this choice of \(\lambda\) and \(\EdgeMultiplier\) is feasible.
\end{proof}

Finite bounds still never reach \(\tfrac12\sqrt{\Biclique}\).
For \(Y\) supported on \(E\), the quadratic form of \(\frac12\sqrt \Biclique\,\Identity-\CostMatrix+Y\) vanishes at the vector \(\frac1{\sqrt2}\left(e_{A}+e_{B}\right)\) of every maximum pair \(\left(A,B\right)\), whereas the form of a matrix of \(\PolyaCone[n]{r}\) with positive diagonal is positive at every nonzero nonnegative vector.

\begin{theorem}[]\protect\hypertarget{thm-polya-inexact}{}\label{thm-polya-inexact}

If \(\Biclique\left(G\right)\ge1\), then \(\Polya{r}\left(G\right)>\tfrac12\sqrt{\Biclique\left(G\right)}\) for every \(r\).

\end{theorem}

\begin{proof}
If \hyperref[eq-zeta]{(\ref*{eq-zeta})} is infeasible, there is nothing to prove.
Otherwise it is a linear program bounded below by \(\frac12\sqrt \Biclique\); hence its minimum is attained, and it suffices to show that every feasible point has \(\lambda>\frac12\sqrt \Biclique\).
Let \(M=\lambda \Identity-\CostMatrix+Y\in\PolyaCone[n]{r}\) with \(Y\) supported on \(E\), so that \(\lambda\ge\frac12\sqrt \Biclique>0\).
The polynomial \(f=\left(\sum_{i}\FormVariable_{i}^{2}\right)^{r}p_{M}\) has only the monomials \(\FormVariable^{2\nu}\), with nonnegative coefficients, and its coefficient of \(\FormVariable_{i}^{2r+4}\) is \(M_{ii}=\lambda\); thus \(f\left(\FormVariable\right)\ge\lambda\sum_{i}\FormVariable_{i}^{2r+4}>0\) for every nonzero \(\FormVariable\).
Let \(u=\frac1{\sqrt2}\left(e_{A}+e_{B}\right)\) for a maximum pair \(\left(A,B\right)\), so that \(u^{\Transpose}\CostMatrix u=\frac12\sqrt \Biclique\), \(\|u\|=1\), and \(u^{\Transpose}Yu=0\) because no edge joins two vertices of \(A\cup B\).
With \(\FormVariable_{i}=\sqrt{u_{i}}\), we get \(0<f\left(\FormVariable\right)=\left(\sum_{i}u_{i}\right)^{r}u^{\Transpose}Mu=\left(\sum_{i}u_{i}\right)^{r}\left(\lambda-\frac12\sqrt \Biclique\right)\); hence \(\lambda>\frac12\sqrt \Biclique\).
\end{proof}

The counterpart of \hyperref[thm-polya-inexact]{Theorem~\ref*{thm-polya-inexact}} for the stability number is the result of Peña et al. \citep{pena2007computing} that the Pólya bounds on \(\Independence\left(G\right)\) are never exact.
As \(\left\langle \CostMatrix,X\right\rangle\) couples only the two sides, the linear programs become feasible at a level set by \(\Biindependence\left(G\right)\), which counts only pairs with both sides nonempty, instead of \(\Independence\left(G\right)\).
Since \(\Biclique\) is an integer and \(\Polya{r}\) decreases to \(\tfrac12\sqrt \Biclique\), \(\Biclique=\lfloor4\left(\Polya{r}\right)^{2}\rfloor\) for all large \(r\), just as rounding down the Pólya bound on \(\Independence\) gives \(\Independence\) at every level \(r\ge\Independence^{2}\) \citep[Theorem 4.1]{deklerk2002approximation}.
However, we have no explicit level such as \(\Independence^{2}\), and no closed form for \(\Polya{r}\) like the one of Peña et al.~for \(\Independence\).

\begin{refremark}
The inexactness in \hyperref[thm-polya-inexact]{Theorem~\ref*{thm-polya-inexact}} is specific to the Pólya cones.
The bound \(\PolyaSOS{1}\) of the level-1 cone \(\PolyaSOSCone[n]{1}\) of Parrilo \citep{parrilo2000structured} is exact on \(\Value{rhooneexact}\) of the \(\Value{graphs}\) random graphs of \hyperref[sec-bounds]{Section~\ref*{sec-bounds}}.
The larger cones \(\PolyaDSOSCone[n]{r}\) built from the diagonally dominant sums of squares (DSOS) of Ahmadi and Majumdar \citep{ahmadi2019dsos} give finite bounds from level \(0\) on, and these bounds can be exact.
At level \(0\) the bound equals half the maximum degree of \(\Complement{G}\).
At that level, however, neither these cones nor the cones \(\PolyaSDSOSCone[n]{0}\) built from the scaled diagonally dominant sums of squares (SDSOS) give a smaller penalty than the spectral threshold (\ref{sec-dsos}).

\label{rem-polya}

\end{refremark}

\section{Computational Results}\label{sec-numerics}

The experiments address three questions: how much smaller our threshold is than the earlier ones (\hyperref[sec-exp-threshold]{Section~\ref*{sec-exp-threshold}}), whether the smaller penalty helps a local ascent find a maximum pair (\hyperref[sec-ascent]{Section~\ref*{sec-ascent}}), and how close the bounds of \hyperref[sec-relax]{Section~\ref*{sec-relax}} come to \(\tfrac12\sqrt{\Biclique}\) (\hyperref[sec-bounds]{Section~\ref*{sec-bounds}}, \hyperref[sec-structured]{Section~\ref*{sec-structured}}).
A caveat is that rounding the bound \(\PolyaSOS{1}\) of \hyperref[sec-rho]{Section~\ref*{sec-rho}} already recovers \(\Biclique\) on every random graph of \hyperref[sec-bounds]{Section~\ref*{sec-bounds}}.
The penalty of \hyperref[thm-penalty]{Theorem~\ref*{thm-penalty}} is therefore not a cheaper way to compute \(\Biclique\), and its role is to find a maximum pair, which the rounded bound then certifies.

All graphs are small, with at most \(40\) vertices, so that \(\Biclique\) can be computed by enumerating the subsets of \(V_{1}\).
The programs were modeled in \texttt{CVXPY} \citep{diamond2016cvxpy} and solved with the interior-point solver \texttt{MOSEK} \citep{mosek2025python}.
We call a bound exact if it differs from \(\tfrac12\sqrt \Biclique\) by less than \(\sci{\Value{exacttol}}\), and we report gaps in relative terms.
The computed absolute gaps of \(\PolyaSOS{1}\) are either below \(\sciatmost{2}{\Value{gapbelow}}\) or above \(\sciatleast{2}{\Value{gapabove}}\), so that its count of exact graphs does not depend on this tolerance.
The code is available at \href{https://github.com/eminyous/bipen/tree/v0.1.4}{\texttt{github.com/eminyous/bipen}}, version 0.1.4; it was written with the assistance of Claude (Anthropic) and checked by the authors.

\subsection{Penalty Thresholds}\label{sec-exp-threshold}

We compare the smallest exact penalty \(\SmallestPenalty\) with four thresholds, named after the quantities they are built on: \(\PenaltyThreshold{g}=\frac12\left(1+\sqrt \Biclique\right)\) of \hyperref[thm-penalty]{Theorem~\ref*{thm-penalty}}, \(\PenaltyThreshold{m}=\frac12\left(1+\max\left(m_{1},m_{2}\right)\right)\) and \(\PenaltyThreshold{E}=\frac12\left(1+\sqrt{\left|E\left(\Complement{G}\right)\right|}\right)\) of Gillis and Glineur \citep{gillis2014continuous}, and the spectral threshold \(\PenaltyThreshold{s}=\frac12\left(1+\|\OnesMatrix-\Biadjacency\|_{2}\right)\).
They satisfy \(\PenaltyThreshold{g}\le \PenaltyThreshold{s}\le \PenaltyThreshold{E}\le \PenaltyThreshold{m}\) by \hyperref[sec-local]{Section~\ref*{sec-local}}.
We drew \(100\) random graphs, with every edge present independently with probability \(p\), for each size \(2+3\), \(3+3\), \(3+4\), \(4+4\), \(4+5\) and \(5+5\) and each \(p=0.2\), \(0.4\), \(0.6\) and \(0.8\).
Discarding the empty and the complete ones left \(\Value{penaltygraphs}\) graphs, for which we computed \(\SmallestPenalty\left(G\right)\) by bisection, evaluating \hyperref[eq-mu]{(\ref*{eq-mu})} exactly as described in \ref{sec-comp}.
The computed \(\SmallestPenalty\left(G\right)\) never exceeded \(\PenaltyThreshold{g}\), as \hyperref[thm-penalty]{Theorem~\ref*{thm-penalty}} requires, and was at most \(\atmost{3}{\Value{penaltyratio}}\) times \(\PenaltyThreshold{g}\) (\hyperref[tbl-threshold]{Table~\ref*{tbl-threshold}}).
The threshold \(\PenaltyThreshold{g}\) is well below \(\PenaltyThreshold{m}\) at every size, and below \(\PenaltyThreshold{E}\) by a margin that grows with the size of the graph, while the spectral threshold \(\PenaltyThreshold{s}\) lies between \(\PenaltyThreshold{g}\) and \(\PenaltyThreshold{E}\).
For the graphs \(N_{k}\) of \hyperref[prp-sharp]{Proposition~\ref*{prp-sharp}} with \(2\le k\le5\), the computed values are within \(\sciatmost{1}{\Value{nearbicliqueerror}}\) of \(\SmallestPenalty\left(N_{k}\right)=\frac k2\).

\begin{table}[tbp]
\centering\small
\caption{Ratios of the smallest exact penalty $\SmallestPenalty$ and of the thresholds $\PenaltyThreshold{g}=\frac12\left(1+\sqrt \Biclique\right)$, $\PenaltyThreshold{m}=\frac12\left(1+\max\left(m_{1},m_{2}\right)\right)$, $\PenaltyThreshold{E}=\frac12\left(1+\sqrt{\left|E\left(\Complement{G}\right)\right|}\right)$ and $\PenaltyThreshold{s}=\frac12\left(1+\|\OnesMatrix-\Biadjacency\|_{2}\right)$, computed graph by graph, with their mean and maximum or minimum over the \# random graphs of each size and their value on the graphs $N_{k}$.}\label{tbl-threshold}
\readrows\randomrows{filter strcmp={\colfamily}{random}}{\TableData{thresholds}}{$\colgraphs$ & $\colcount$ & $\fixed{3}{\colgammastaroveronemean}$ & $\fixed{3}{\colgammastaroveronemax}$ & $\fixed{3}{\colgammaoneovertwomean}$ & $\fixed{3}{\colgammaoneovertwomin}$ & $\fixed{3}{\colgammaoneoverthreemean}$ & $\fixed{3}{\colgammaoneoverthreemin}$ & $\fixed{3}{\colgammaoneoverfourmean}$ & $\fixed{3}{\colgammaoneoverfourmin}$}
\readrows\nkrows{filter strcmp={\colfamily}{nearbiclique}}{\TableData{thresholds}}{$N_{\colgraphs}$ & & \csvexpnot\multicolumn{2}{c}{$\fixed{3}{\colgammastaroveronemean}$} & \csvexpnot\multicolumn{2}{c}{$\fixed{3}{\colgammaoneovertwomean}$} & \csvexpnot\multicolumn{2}{c}{$\fixed{3}{\colgammaoneoverthreemean}$} & \csvexpnot\multicolumn{2}{c}{$\fixed{3}{\colgammaoneoverfourmean}$}}

\begin{tabular}{@{}lrrrrrrrrr@{}}
\toprule
& & \multicolumn{2}{c}{$\SmallestPenalty/\PenaltyThreshold{g}$} & \multicolumn{2}{c}{$\PenaltyThreshold{g}/\PenaltyThreshold{m}$} & \multicolumn{2}{c}{$\PenaltyThreshold{g}/\PenaltyThreshold{E}$} & \multicolumn{2}{c}{$\PenaltyThreshold{g}/\PenaltyThreshold{s}$}\\
\cmidrule(lr){3-4}\cmidrule(lr){5-6}\cmidrule(lr){7-8}\cmidrule(l){9-10}
$G$ & \# & mean & max & mean & min & mean & min & mean & min\\
\midrule
\randomrows
\midrule
\nkrows
\bottomrule
\end{tabular}
\end{table}

\subsection{Local Ascent}\label{sec-ascent}

We ran block-coordinate ascent on \hyperref[eq-sigma]{(\ref*{eq-sigma})}, which maximizes exactly over one side at a time.
Given \(z\), the update is \(y\leftarrow\|\left(\PenaltyMatrix{\gamma} z\right)_{+}\|^{-1}\left(\PenaltyMatrix{\gamma} z\right)_{+}\), the normalized positive part of \(\PenaltyMatrix{\gamma} z\).
When this positive part is zero, the update is \(y\leftarrow e_{i}\) for a vertex \(i\) maximizing \(\left(\PenaltyMatrix{\gamma} z\right)_{i}\), and we call this second rule the fallback.
The updates of \(z\) are the same with the two sides exchanged.
We ran the ascent on the first \(20\) random graphs of \hyperref[sec-bounds]{Section~\ref*{sec-bounds}} of each size \(10+10\), \(15+15\) and \(20+20\) and each edge probability \(p=0.3\), \(0.5\) and \(0.7\), from two kinds of starts.
The random starts are \(100\) points with entries uniform in \(\left[0,1\right]\), drawn with \texttt{numpy.random.default\_rng(s)} for the graph of seed index \(s\).
The vertex starts are the \(n\) vertices, where the start \(z=e_{j}\) makes the first \(y\) uniform on the non-neighbors of \(j\), and likewise for \(y=e_{i}\).
The penalties are the four thresholds of \hyperref[sec-exp-threshold]{Section~\ref*{sec-exp-threshold}}, where \(\PenaltyThreshold{g}\) is available by the rounding of \hyperref[sec-bounds]{Section~\ref*{sec-bounds}}.
A run succeeds when its final supports form a maximum pair, which for \(\gamma\ge \PenaltyThreshold{g}\) means that it ends at the value \(\sqrt \Biclique\).
For every size, edge probability, kind of start and penalty, the ascent stops after at most \(\atmost{0}{\Value{ascentsteps}}\) rounds of the two updates on average, and it therefore mostly maps the starting point to a nearby pair.

In every row of \hyperref[tbl-ascent]{Table~\ref*{tbl-ascent}}, the ascent with \(\PenaltyThreshold{g}\) succeeds from a larger fraction of the starts than with the three larger penalties, and on at least as many graphs.
From random starts at \(p=0.3\) and \(0.5\) the three larger penalties never succeed beyond \(10+10\), while \(\PenaltyThreshold{g}\) still finds a maximum pair from some start on \(\Value{ascentgraphs}\) of the \(20\) graphs with \(15+15\) vertices and \(p=0.3\).
The spectral threshold \(\PenaltyThreshold{s}\), the smallest threshold that needs neither \(\Biclique\) nor a semidefinite program, helps little, as beyond \(10+10\) and for \(p\le0.5\) it succeeds on at most one graph for each size and edge probability, from either kind of start.

The random starts at \(p=0.7\) are left out of \hyperref[tbl-ascent]{Table~\ref*{tbl-ascent}} because they mostly measure the fallback.
Their first step uses the fallback in at least \(\atleast{0}{\Value{densefallback}}\%\) of the runs for every penalty, and the vertex \(i\) that it selects, which maximizes \(\left(\PenaltyMatrix{\gamma} z\right)_{i}=\OnesVector^{\Transpose}z-2\gamma\left(\Biadjacency z\right)_{i}\), does not depend on \(\gamma\).
On these runs the three larger penalties give identical results, and \(\PenaltyThreshold{g}\) succeeds from a somewhat larger fraction of the starts at each size.
From the vertex starts the fallback occurs in at most \(\atmost{1}{\Value{vertexfallback}}\%\) of the runs.

With the four penalties, the unsuccessful runs end at a maximal pair that is not maximum, except for at most \(\atmost{1}{\Value{edgeend}}\%\) of the runs of any size and edge probability, which end with an edge of \(G\) in the support.
In all, \(\Value{edgeruns}\) runs end in this way, at fixed points of the block updates, and by \hyperref[prp-local]{Proposition~\ref*{prp-local}} the \(\Value{edgerunspositive}\) of them with positive value are not local maximizers.
No run ends at a biindependent pair that is not maximal, because at a fixed point with biindependent supports every vertex of \(V_{1}\) without a neighbor in \(\nz\left(z\right)\) has the largest entry \(\OnesVector^{\Transpose}z\) of \(\PenaltyMatrix{\gamma} z\) and hence lies in \(\nz\left(y\right)\), and likewise for \(V_{2}\).
A greedy repair of the final supports, which removes the vertex with the most neighbors on the other side until the pair is biindependent and then adds vertices of smallest degree while it stays biindependent, changes no result with these penalties.

A smaller penalty helps because each step keeps in \(y\) only the vertices \(i\) with \(2\gamma\left(\Biadjacency z\right)_{i}<\OnesVector^{\Transpose}z\).
A larger penalty therefore discards more vertices at once, and the ascent settles at a maximal pair near a small support.
For \(10+10\) and \(p=0.3\), the first step from a random start keeps on average \(\fixed{0}{\Value{supportone}}\%\) of \(V_{1}\) with \(\PenaltyThreshold{g}\) and \(\fixed{0}{\Value{supporttwo}}\%\) with \(\PenaltyThreshold{m}\), and it uses the fallback in \(\fixed{0}{\Value{fallbackone}}\%\) and \(\fixed{0}{\Value{fallbacktwo}}\%\) of the runs.
Penalties below \(\PenaltyThreshold{g}\) keep even more vertices but lose the guarantee of \hyperref[prp-local]{Proposition~\ref*{prp-local}}.
With some of them the ascent succeeds more often (\hyperref[tbl-ascent-below]{Table~\ref*{tbl-ascent-below}} in \ref{sec-comp}), and thus the threshold of \hyperref[thm-penalty]{Theorem~\ref*{thm-penalty}} is not the best penalty for this heuristic.

\begin{table}[tbp]
\centering\small\setlength{\tabcolsep}{5pt}
\caption{Fraction of starts from which the local ascent ends at a maximum pair, and number of graphs (out of $20$) on which some start does, with the penalties $\PenaltyThreshold{g}=\frac12\left(1+\sqrt \Biclique\right)$, $\PenaltyThreshold{s}=\frac12\left(1+\|\OnesMatrix-\Biadjacency\|_{2}\right)$, $\PenaltyThreshold{E}=\frac12\left(1+\sqrt{\left|E\left(\Complement{G}\right)\right|}\right)$ and $\PenaltyThreshold{m}=\frac12\left(1+\max\left(m_{1},m_{2}\right)\right)$.}\label{tbl-ascent}
\readrows\randomrows{filter strcmp={\colstarts}{random}, and filter fp={\colp<0.7}}{\TableData{ascent}}{$\colsize$ & $\colp$ & $\fixed{3}{\colgammaonesuccess}$ & $\colgammaonegraphs$ & $\fixed{3}{\colgammafoursuccess}$ & $\colgammafourgraphs$ & $\fixed{3}{\colgammathreesuccess}$ & $\colgammathreegraphs$ & $\fixed{3}{\colgammatwosuccess}$ & $\colgammatwographs$}
\readrows\vertexrows{filter strcmp={\colstarts}{vertex}}{\TableData{ascent}}{$\colsize$ & $\colp$ & $\fixed{3}{\colgammaonesuccess}$ & $\colgammaonegraphs$ & $\fixed{3}{\colgammafoursuccess}$ & $\colgammafourgraphs$ & $\fixed{3}{\colgammathreesuccess}$ & $\colgammathreegraphs$ & $\fixed{3}{\colgammatwosuccess}$ & $\colgammatwographs$}

\begin{tabular}{@{}lrrrrrrrrr@{}}
\toprule
& & \multicolumn{2}{c}{$\PenaltyThreshold{g}$} & \multicolumn{2}{c}{$\PenaltyThreshold{s}$} & \multicolumn{2}{c}{$\PenaltyThreshold{E}$} & \multicolumn{2}{c}{$\PenaltyThreshold{m}$}\\
\cmidrule(lr){3-4}\cmidrule(lr){5-6}\cmidrule(lr){7-8}\cmidrule(l){9-10}
$G$ & $p$ & starts & graphs & starts & graphs & starts & graphs & starts & graphs\\
\midrule
\multicolumn{10}{@{}l}{\emph{Random starts}}\\
\randomrows
\midrule
\multicolumn{10}{@{}l}{\emph{Vertex starts}}\\
\vertexrows
\bottomrule
\end{tabular}
\end{table}

\subsection{Bounds on Random Graphs}\label{sec-bounds}

We drew \(100\) random bipartite graphs for each edge probability \(p=0.3\), \(0.5\) and \(0.7\) and each of ten sizes \(m_{1}+m_{2}\), and computed the bounds \(\LasserreBicliqueRatio{1}\), \(\ConeBicliqueRatio{\DNN}\) and \(\PolyaSOS{1}\) of \hyperref[sec-relax]{Section~\ref*{sec-relax}} on them.
Six sizes are balanced, \(6+6\), \(8+8\), \(10+10\), \(12+12\), \(15+15\) and \(20+20\), and four are unbalanced, \(5+10\), \(8+16\), \(10+20\) and \(6+24\).
For each size, \hyperref[tbl-random]{Table~\ref*{tbl-random}} gives the number of graphs on which each bound is exact and the relative gaps, and \hyperref[tbl-random-p]{Table~\ref*{tbl-random-p}} in \ref{sec-comp} separates the three edge probabilities.

Adding the constraint \(X\ge0\), that is, passing from \(\LasserreBicliqueRatio{1}\) to \(\ConeBicliqueRatio{\DNN}\), makes the bound exact on more graphs and reduces the mean gap at every size.
The bound \(\PolyaSOS{1}\) is exact on \(\Value{rhooneexact}\) of the \(\Value{graphs}\) graphs and exceeds \(\frac12\sqrt \Biclique\) by at most \(\atmost{2}{\Value{rhoonemaxgap}}\%\) on the other \(\Value{rhooneinexact}\).
The bounds \(\ConeBicliqueRatio{\DNN}\) and \(\PolyaSOS{1}\) are exact least often at \(p=0.3\), where \(\Value{rhooneinexactsparse}\) of the \(\Value{rhooneinexact}\) failures of \(\PolyaSOS{1}\) occur, while the mean gap of \(\ConeBicliqueRatio{\DNN}\) over all sizes is largest at \(p=0.5\).

\begin{table}[tbp]
\centering\small\setlength{\tabcolsep}{3.5pt}
\caption{The bounds $\LasserreBicliqueRatio{1}$, $\ConeBicliqueRatio{\DNN}$ and $\PolyaSOS{1}$ on random bipartite graphs with $|V_{1}|=m_{1}$ and $|V_{2}|=m_{2}$, with the edge probabilities $p=0.3$, $0.5$ and $0.7$ together ({\hyperref[tbl-random-p]{Table~\ref*{tbl-random-p}}} separates them); \# is the number of graphs.
Columns 3--5: number of graphs on which the bound is within $\sci{\Value{exacttol}}$ of $\frac12\sqrt \Biclique$.
Columns 6--11: gap of the bound relative to $\frac12\sqrt \Biclique$, in percent, mean and maximum over the graphs of the row.}\label{tbl-random}
\readrows\balancedrows{filter strcmp={\colblock}{balanced}, and filter strcmp={\colp}{all}}{\TableData{bounds}}{$\colsize$ & \boundscells}
\readrows\unbalancedrows{filter strcmp={\colblock}{unbalanced}, and filter strcmp={\colp}{all}}{\TableData{bounds}}{$\colsize$ & \boundscells}

\begin{tabular}{@{}rrrrrrrrrrr@{}}
\toprule
& & \multicolumn{3}{c}{exact} & \multicolumn{3}{c}{mean gap (\%)} & \multicolumn{3}{c}{max gap (\%)}\\
\cmidrule(lr){3-5}\cmidrule(lr){6-8}\cmidrule(l){9-11}
$m_{1}+m_{2}$ & \boundsheads\\
\midrule
\multicolumn{11}{@{}l}{\emph{balanced, $m_{1}=m_{2}$}}\\
\balancedrows
\midrule
\multicolumn{11}{@{}l}{\emph{unbalanced, $m_{1}<m_{2}$}}\\
\unbalancedrows
\bottomrule
\end{tabular}
\end{table}

\textbf{Rounding.}
Since \(\Biclique\) is an integer, a bound \(b\ge\frac12\sqrt \Biclique\) gives \(\Biclique=\lfloor4b^{2}\rfloor\) as soon as \(b<\frac12\sqrt{\Biclique+1}\), and it therefore need not be exact.
A computed bound falls below \(\frac12\sqrt \Biclique\) by at most \(\sciatmost{1}{\Value{solvererror}}\), and rounding \(4b^{2}+\sci{\Value{roundoffset}}\) down absorbs this error.
Rounded in this way, \(\PolyaSOS{1}\) recovers \(\Biclique\) on \allof{\Value{roundrhoone}}{\Value{graphs}} graphs, the \(\Value{rhooneinexact}\) inexact ones included, whereas \(\ConeBicliqueRatio{\DNN}\) recovers it on \(\Value{roundhonednn}\) graphs and \(\LasserreBicliqueRatio{1}\) on \(\Value{roundhone}\).
As \(\lfloor4b^{2}\rfloor\ge \Biclique\), a pair with \(|A||B|=\lfloor4b^{2}\rfloor\), such as one found by the local ascent of \hyperref[sec-ascent]{Section~\ref*{sec-ascent}}, is maximum.

\textbf{Lasserre bound.}
Laurent et al. \citep[(29)]{laurent2025semidefinite} also define bounds \(\LasserreBiclique{r}\left(G\right)\ge \Biclique\left(G\right)\), and we compare \(\PolyaSOS{1}\) with the second of them through \(\tfrac12\sqrt{\LasserreBiclique{2}\left(G\right)}\).
We do not use the bounds \(\LasserreBicliqueRatio{r}\) of the same paper, since these bound \(\BicliqueRatio\left(G\right)\), which can be smaller than \(\tfrac12\sqrt \Biclique\).
We computed \(\LasserreBiclique{2}\) in moment form, with the moment matrix indexed by the stable sets of size at most \(2\), on the first ten graphs of each balanced size up to \(15+15\) and each edge probability, \(\Value{lasserregraphs}\) graphs in all.
The Lasserre bound was exact on all of them and \(\PolyaSOS{1}\) on \(\Value{lasserrerhooneexact}\), and the Lasserre bound was also faster.
On these small graphs the second Lasserre level is the better bound on \(\Biclique\).
Its cost grows like \(n^{4}\) and is larger on sparse graphs, as its moment matrix involves the stable sets of size at most \(4\), while \(\PolyaSOS{1}\) has \(n\) blocks of order \(n\).
On three graphs with \(20+20\) vertices for each edge probability, nine in all, which give only an indication, the Lasserre bound was slower than \(\PolyaSOS{1}\) for \(p=0.3\), comparable for \(p=0.5\) and still faster for \(p=0.7\).
On one of the graphs with \(p=0.3\) both bounds are inexact and the Lasserre bound is the weaker one, with the value \(\fixedtrim{4}{\Value{weaklasserre}}\) where \(\PolyaSOS{1}\) has the value \(\fixedtrim{4}{\Value{weakrhoone}}\) and \(\tfrac12\sqrt \Biclique=\fixedtrim{4}{\Value{weakhalfsqrtg}}\).
Both bounds are exact on the six graphs with \(p\ge0.5\).

\subsection{Structured Graphs and Pólya Bounds}\label{sec-structured}

Structured graphs show where the relaxations of \hyperref[sec-relax]{Section~\ref*{sec-relax}} are exact and how slowly the Pólya bounds converge.
On the graph of Laurent et al. \citep[Example 1]{laurent2025semidefinite}, which has one edge on \(2+2\) vertices, and on the crowns \(K_{s,s}\) minus a perfect matching, \(5\le s\le12\), all bounds equal \(\frac12\sqrt \Biclique\), namely \(\frac12\sqrt2\) and \(\frac12\) respectively, consistent with \(\LasserreBicliqueRatio{1}=\frac12\sqrt \Biclique\) on these graphs \citep[Example 1 and Section 5]{laurent2025semidefinite}.
On the perfect matchings \(sK_{2}\), \(2\le s\le12\), all bounds equal \(\frac12\sqrt \Biclique=\frac s4\) for even \(s\).
For odd \(s\) (\hyperref[tbl-matchings]{Table~\ref*{tbl-matchings}}), \(\LasserreBicliqueRatio{1}\) and \(\ConeBicliqueRatio{\DNN}\) equal \(\frac s4>\frac12\sqrt \Biclique\), as in \hyperref[prp-matching]{Proposition~\ref*{prp-matching}}, and \(\PolyaSOS{1}\) is exact only for \(s=3\).

\begin{table}[tbp]
\centering\small
\caption{The bounds $\LasserreBicliqueRatio{1}$, $\ConeBicliqueRatio{\DNN}$ and $\PolyaSOS{1}$ on the perfect matchings $sK_{2}$, $2\le s\le12$, on which one of them is not exact, that is, for odd $s$; there $\frac12\sqrt \Biclique=\frac14\sqrt{s^{2}-1}$, while for even $s$ all bounds equal $\frac s4$.}\label{tbl-matchings}
\readrows\matchingrows{filter fp={max(\colhone,\colhonednn,\colrhoone)-\colhalfsqrtg>=\Value{exacttol}}}{\TableData{matchings}}{$\colr$ & $\fixed{4}{\colhalfsqrtg}$ & $\fixed{4}{\colhone}$ & $\fixed{4}{\colhonednn}$ & $\fixed{4}{\colrhoone}$}

\begin{tabular}{@{}rrrrr@{}}
\toprule
$s$ & $\frac12\sqrt \Biclique$ & $\LasserreBicliqueRatio{1}$ & $\ConeBicliqueRatio{\DNN}$ & $\PolyaSOS{1}$\\
\midrule
\matchingrows
\bottomrule
\end{tabular}
\end{table}

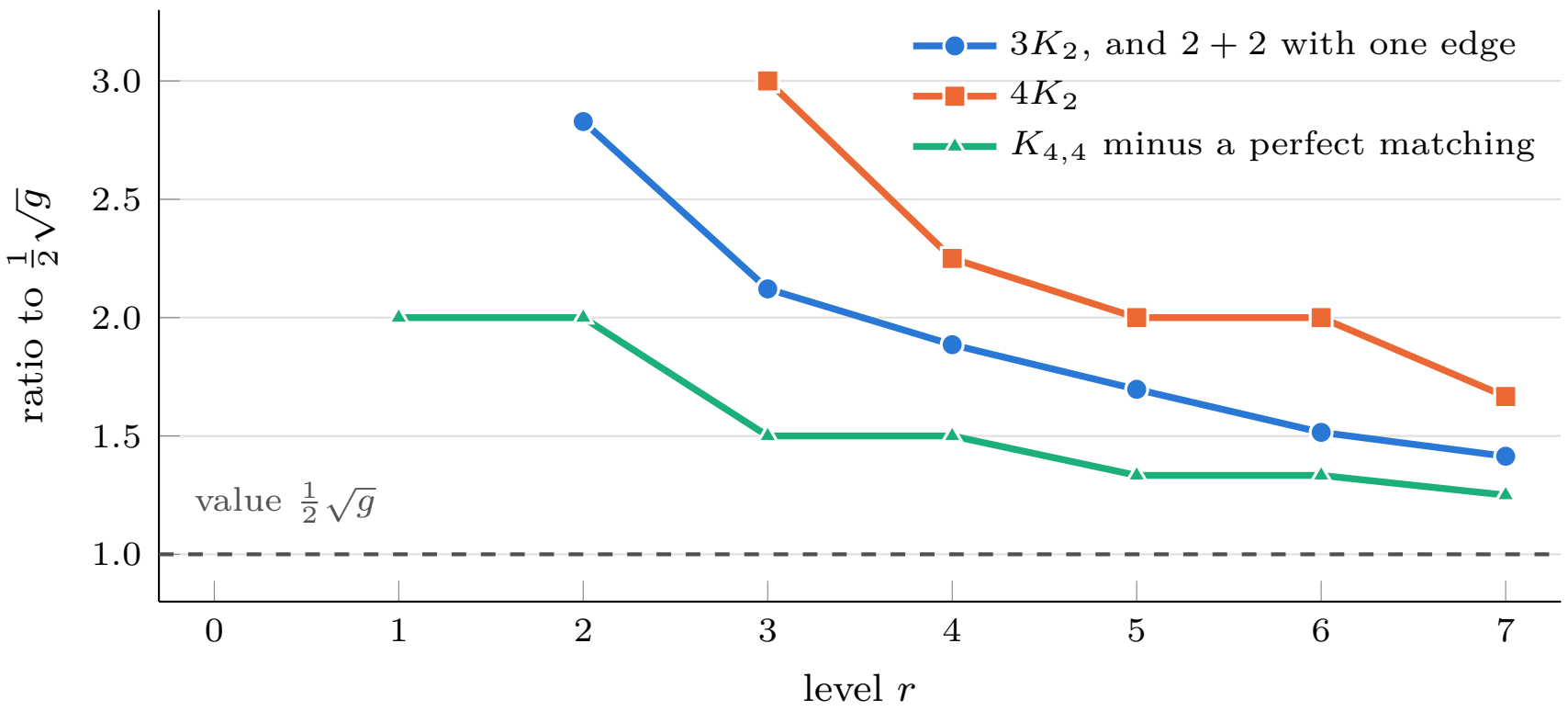
\begin{figure}[tbp]
\centering
\pgfplotstableread[col sep=comma]{\FigureData{polya}}\polyatable
\begin{tikzpicture}
\begin{axis}[width=4.6in, height=2.3in, axis lines*=left, xmin=-0.3, xmax=7.3, ymin=0.8, ymax=3.3, xtick={0,...,7}, ytick={1,1.5,...,3},
  yticklabel style={/pgf/number format/.cd, fixed, fixed zerofill, precision=1},
  xlabel={level $r$}, ylabel={ratio to $\frac12\sqrt \Biclique$}, label style={font=\footnotesize}, tick label style={font=\scriptsize},
  ymajorgrids, grid style={black!12, line width=0.4pt}, legend cell align=left,
  legend style={draw=none, fill=none, font=\scriptsize, at={(1,1)}, anchor=north east},
  every axis plot/.append style={line width=1.4pt, mark size=2.25pt}]
\addplot[plotmuted, dashed, line width=0.8pt, mark=none, domain=-0.3:7.3, samples=2, forget plot] {1};
\node[anchor=south west, font=\scriptsize, text=plotmuted] at (axis cs:-0.2,1.06) {value $\frac12\sqrt \Biclique$};
\addplot[plotblue, mark=*, mark options={draw=white, line width=0.6pt, fill=plotblue}] table[x=level, y=matchingthree] {\polyatable};
\addlegendentry{$3K_{2}$, and $2+2$ with one edge}
\addplot[plotorange, mark=square*, mark options={draw=white, line width=0.6pt, fill=plotorange}] table[x=level, y=matchingfour] {\polyatable};
\addlegendentry{$4K_{2}$}
\addplot[plotgreen, mark=triangle*, mark options={draw=white, line width=0.6pt, fill=plotgreen}] table[x=level, y=crownfour] {\polyatable};
\addlegendentry{$K_{4,4}$ minus a perfect matching}
\end{axis}
\end{tikzpicture}
\caption{The Pólya linear programming bounds $\Polya{r}$ divided by $\tfrac12\sqrt \Biclique$, by level $r$, on four small graphs.
Each program is infeasible below its first marker, at level $\Biindependence-1$, where $\Biindependence$ is the largest $|A|+|B|$ over the biindependent pairs with both sides nonempty}\label{fig-polya}
\end{figure}

\textbf{Pólya bounds.}
The bounds of \hyperref[eq-zeta]{(\ref*{eq-zeta})} approach \(\tfrac12\sqrt \Biclique\) slowly on the four small graphs of \hyperref[fig-polya]{Figure~\ref*{fig-polya}}.
The curves of \(3K_{2}\) and of its induced subgraph with one edge on \(2+2\) vertices coincide numerically.
On \(3K_{2}\), the bound at level \(\Value{polyalevel}\) is still \(\fixed{2}{\Value{polyaratio}}\) times \(\tfrac12\sqrt \Biclique\), while \(\PolyaSOS{1}\left(3K_{2}\right)\) is already exact.

\section{Conclusions}\label{sec-conclusion}

The penalty that makes the continuous formulation of the maximum-edge biclique problem exact is set by the square root of the optimal value, not by the size of the graph.
Our threshold cannot be lowered by more than one half in general.
The same penalty is an optimal dual multiplier in an exact pair of copositive and completely positive programs, whose relaxations make it computable.

The computational evidence is limited to graphs with at most forty vertices.
On these graphs a rounded semidefinite bound already gave the optimal value, and the penalty served only to find a maximum pair.
From random starting points on the graphs with thirty and forty vertices and edge probability at most one half, the ascent reached maximum pairs with this penalty and never with the larger thresholds, but only from a small fraction of the starting points.
Some penalties below the threshold, which carry no guarantee, did better.

Three questions remain open.
From some level on, rounding the Pólya bounds recovers the optimal value, but we know neither an explicit such level nor a closed form for these bounds.
The smallest exact penalty can equal half the square root of the optimal value, but we computed it only on graphs with at most eleven vertices, and its behavior on larger graphs is unknown.
How to choose the penalty for the local ascent, given that values below the threshold can do better, is open as well.

\section*{Statements and Declarations}\label{statements-and-declarations}
\addcontentsline{toc}{section}{Statements and Declarations}

\subsection*{Funding}\label{funding}
\addcontentsline{toc}{subsection}{Funding}

This work was supported by Bissan Ghaddar's Natural Sciences and Engineering Research Council of Canada Discovery Grant RGPIN-2025-04585 and John Thompson Chair Fellowship.

\subsection*{Competing Interests}\label{competing-interests}
\addcontentsline{toc}{subsection}{Competing Interests}

The authors declare that they have no known competing financial interests or personal relationships that could have appeared to influence the work reported in this paper.

\subsection*{Author Contributions}\label{author-contributions}
\addcontentsline{toc}{subsection}{Author Contributions}

Youssouf Emine: conceptualization, methodology, formal analysis, software, writing of the original draft.
Bissan Ghaddar: supervision, funding acquisition, review and editing.
Both authors read and approved the final manuscript.

\subsection*{Data and Code Availability}\label{data-and-code-availability}
\addcontentsline{toc}{subsection}{Data and Code Availability}

The numerical results of this paper are generated by the code available at \href{https://github.com/eminyous/bipen}{\texttt{github.com/eminyous/bipen}} under the MIT license, in the version given in \hyperref[sec-numerics]{Section~\ref*{sec-numerics}}; no other data were used.

\subsection*{Use of Generative AI}\label{use-of-generative-ai}
\addcontentsline{toc}{subsection}{Use of Generative AI}

During the preparation of this work, the authors used Claude (Anthropic) to assist with writing the code for the numerical experiments, correcting the English, clarifying the wording, and searching the literature.
All theorems and proofs are the authors' own work.
After using this tool, the authors reviewed and edited the content as needed and take full responsibility for the content of the published article.

\appendix
\appendixnumbering

\section{Bipartite Graphs on Five Vertices}\label{sec-five}

This appendix proves \hyperref[prp-small]{Proposition~\ref*{prp-small}} for \(n=5\).
Let \(\GenericCone\) be one of the cones \(\PSD[n]\), \(\DNN[n]\) and \(\CP[n]\); each is closed under taking principal submatrices and under bordering by a zero row and column.
A vertex is \emph{full} if it is adjacent to every vertex of the other side.

\begin{lemma}[]\protect\hypertarget{lem-full}{}\label{lem-full}

Let \(v\in V_{1}\) be full, with \(|V_{1}|\ge2\).
Then \(\Biclique\left(G-v\right)=\Biclique\left(G\right)\) and \(\ConeBicliqueRatio{\GenericCone}\left(G-v\right)=\ConeBicliqueRatio{\GenericCone}\left(G\right)\).
The same holds for a full vertex \(v\in V_{2}\) with \(|V_{2}|\ge2\).

\end{lemma}

\begin{proof}
A biindependent pair \(\left(A,B\right)\) of \(G\) with \(v\in A\) has \(B=\varnothing\), and thus \(|A||B|=0\), and the other biindependent pairs of \(G\) are those of \(G-v\); hence \(\Biclique\left(G-v\right)=\Biclique\left(G\right)\).
Bordering a matrix that is feasible for \hyperref[eq-h1K]{(\ref*{eq-h1K})} on \(G-v\) by a zero row and column for \(v\) gives a matrix feasible on \(G\) with the same value; hence \(\ConeBicliqueRatio{\GenericCone}\left(G\right)\ge \ConeBicliqueRatio{\GenericCone}\left(G-v\right)\).
Conversely, let \(X\) be feasible on \(G\) with \(\left\langle \CostMatrix,X\right\rangle>0\), and let \(X^{\prime}\) be its principal submatrix on \(V\setminus\left\{v\right\}\).
Since \(v\) is adjacent to every vertex of \(V_{2}\), the entries \(X_{vj}\), \(j\in V_{2}\), vanish; thus \(\left\langle \CostMatrix,X\right\rangle=\left\langle \CostMatrix^{\prime},X^{\prime}\right\rangle\), where \(\CostMatrix^{\prime}\) is the matrix \(\CostMatrix\) of \(G-v\).
Hence \(X^{\prime}\ne0\); thus \(t=\tr\left(X^{\prime}\right)\) lies in \(\left(0,1\right]\), and \(\frac1tX^{\prime}\) is feasible on \(G-v\) with value \(\frac1t\left\langle \CostMatrix,X\right\rangle\ge\left\langle \CostMatrix,X\right\rangle\).
As \(\ConeBicliqueRatio{\GenericCone}\left(G-v\right)\ge0\), this gives \(\ConeBicliqueRatio{\GenericCone}\left(G-v\right)\ge \ConeBicliqueRatio{\GenericCone}\left(G\right)\).
\end{proof}

\begin{proof}[Proof of Proposition \ref{prp-small} for \(n=5\)]
Delete full vertices one at a time, as long as the side of the deleted vertex keeps another vertex; by \hyperref[lem-full]{Lemma~\ref*{lem-full}} this changes neither \(\Biclique\) nor \(\ConeBicliqueRatio{\DNN}\left(G\right)\).
If fewer than five vertices remain, the case \(n\le4\) applies.
Otherwise no side with at least two vertices contains a full vertex.
If \(|V_{1}|=1\), a neighbor of the vertex of \(V_{1}\) would be a full vertex of \(V_{2}\); thus \(E=\varnothing\), and likewise if \(|V_{2}|=1\).
If \(E=\varnothing\), then \(\ConeBicliqueRatio{\DNN}\left(G\right)\le \LasserreBicliqueRatio{1}\left(G\right)\le\lambda_{\max}\left(\CostMatrix\right)=\frac12\sqrt{m_{1}m_{2}}=\frac12\sqrt \Biclique\), because \(\left\langle \CostMatrix,X\right\rangle\le\lambda_{\max}\left(\CostMatrix\right)\tr\left(X\right)\) for every \(X\succeq0\).
There remains the case \(|V_{1}|=2\) and \(|V_{2}|=3\), up to exchanging the sides, with \(V_{1}=\left\{v_{11},v_{12}\right\}\) and \(V_{2}=\left\{v_{21},v_{22},v_{23}\right\}\).
No vertex of \(V_{2}\) is adjacent to both \(v_{11}\) and \(v_{12}\), and neither of them is adjacent to every vertex of \(V_{2}\); hence the neighborhoods \(N\left(v_{11}\right)\) and \(N\left(v_{12}\right)\) are disjoint and have at most two vertices each.
Up to symmetry, \(\left(|N\left(v_{11}\right)|,|N\left(v_{12}\right)|\right)\) is \(\left(0,0\right)\), a graph without edges, or one of \(\left(0,1\right)\), \(\left(0,2\right)\), \(\left(1,1\right)\) and \(\left(1,2\right)\), the graphs \(G_{1},\dots,G_{4}\) of \hyperref[tab-five]{Table~\ref*{tab-five}}.

For each \(G_{i}\), \hyperref[tab-five]{Table~\ref*{tab-five}} gives \(\lambda=\frac12\sqrt \Biclique\), a number \(\EdgeMultiplier\) and a symmetric nonnegative matrix \(Z\).
Let \(Y\) be the symmetric matrix with entries \(\EdgeMultiplier\) on the edges and \(0\) elsewhere.
Then \(S=\lambda \Identity-\CostMatrix+Y-Z\) is positive semidefinite, as its eigenvalues show.
If \(X\in\DNN[5]\) is feasible for \hyperref[eq-h1K]{(\ref*{eq-h1K})}, then \(\left\langle Y,X\right\rangle=0\) because \(X_{E}=0\); thus
\[
\left\langle \CostMatrix,X\right\rangle=\lambda\tr\left(X\right)-\left\langle S,X\right\rangle-\left\langle Z,X\right\rangle\le\lambda,
\]
since \(S\) and \(X\) are positive semidefinite and \(Z\) and \(X\) are nonnegative.
Hence \(\ConeBicliqueRatio{\DNN}\left(G_{i}\right)\le\frac12\sqrt \Biclique\), and equality holds by \hyperref[thm-cp]{Theorem~\ref*{thm-cp}}.
For \(G_{1}\), \(G_{2}\) and \(G_{4}\) the matrix \(Z\) is zero, and thus the certificate even shows \(\LasserreBicliqueRatio{1}\left(G_{i}\right)=\frac12\sqrt \Biclique\); only \(G_{3}\) needs the constraint \(X\ge0\).
\end{proof}

\begin{table}[H]
\centering\footnotesize\setlength{\tabcolsep}{3.5pt}\renewcommand{\arraystretch}{1.6}
\caption{Certificates for $G_{1},\dots,G_{4}$; the matrix $Z$ is symmetric, with all other entries $0$.}\label{tab-five}
\begin{tabular}{@{}lllllll@{}}
\toprule
Graph & Edges & $\Biclique$ & $\lambda$ & $\EdgeMultiplier$ & $Z$ & Eigenvalues of $S$\\
\midrule
$G_{1}$ & $v_{12}v_{23}$ & $4$ & $1$ & $1$ & $0$ & $0$, $1$, $2$, $1\pm\frac12\sqrt2$\\
$G_{2}$ & $v_{12}v_{22}$, $v_{12}v_{23}$ & $3$ & $\frac12\sqrt3$ & $\frac34$ & $0$ & $0$, $\frac12\sqrt3$, $\sqrt3$, $\frac14\sqrt3\left(2\pm\sqrt2\right)$\\
$G_{3}$ & $v_{11}v_{23}$, $v_{12}v_{22}$ & $2$ & $\frac12\sqrt2$ & $1$ & $Z_{v_{22}v_{23}}=\frac12\sqrt2$ & $0$, $0$, $0$, $\sqrt2$, $\frac32\sqrt2$\\
$G_{4}$ & $v_{11}v_{23}$, $v_{12}v_{21}$, $v_{12}v_{22}$ & $2$ & $\frac12\sqrt2$ & $\frac12$ & $0$ & $0$, $\frac12\sqrt2$, $\sqrt2$, $\frac12\left(\sqrt2\pm1\right)$\\
\bottomrule
\end{tabular}
\end{table}

\section{DSOS and SDSOS Approximations of the Copositive Dual}\label{sec-dsos}

This appendix replaces \(\PolyaCone[n]{r}\) in \hyperref[eq-zeta]{(\ref*{eq-zeta})} by the inner approximations of the copositive cone of \citep{ahmadi2019dsos}, which restrict the Gram matrices of sums of squares.

A symmetric matrix \(\Gamma\) with rows and columns indexed by a finite set is \emph{diagonally dominant} if \(\Gamma_{ii}\ge\sum_{j\ne i}|\Gamma_{ij}|\) for every index \(i\).
It is a \emph{sum of \(2\times2\) blocks} if it is a sum of positive semidefinite matrices each of which vanishes outside the rows and columns of at most two indices.
Nonnegative diagonal matrices have both properties, and both are preserved under sums, nonnegative multiples, principal submatrices, bordering by zero rows and columns, and relabeling of the indices.
A sum of \(2\times2\) blocks is positive semidefinite, and a diagonally dominant matrix \(\Gamma\) is a sum of \(2\times2\) blocks, as it is the nonnegative combination
\[
\Gamma=\sum_{i<j}|\Gamma_{ij}|\left(e_{i}+\sgn\left(\Gamma_{ij}\right)e_{j}\right)\left(e_{i}+\sgn\left(\Gamma_{ij}\right)e_{j}\right)^{\Transpose}+\sum_{i}\left(\Gamma_{ii}-\sum_{j\ne i}|\Gamma_{ij}|\right)e_{i}e_{i}^{\Transpose},
\]
where \(\sgn\left(\Gamma_{ij}\right)\) is the sign of \(\Gamma_{ij}\).

For \(d\ge1\) let \(\left[\FormVariable\right]_{d}\) be the vector of the monomials \(\FormVariable^{\nu}\), \(\nu\in\Naturals^{n}_{d}\).
A form \(q\) of degree \(2d\) in \(\FormVariable\in\Reals^{n}\) is DSOS (a diagonally dominant sum of squares) if it is a nonnegative combination of squares \(\left(\FormVariable^{\nu}\right)^{2}\) and \(\left(\FormVariable^{\nu}\pm \FormVariable^{\nu^{\prime}}\right)^{2}\), and SDSOS (a scaled diagonally dominant sum of squares) if it is a sum of squares \(\left(a \FormVariable^{\nu}+b\FormVariable^{\nu^{\prime}}\right)^{2}\) with \(a,b\in\Reals\), where \(\nu,\nu^{\prime}\in\Naturals^{n}_{d}\) \citep{ahmadi2019dsos}.
Monomials of other degrees cannot occur in such a decomposition, as the parts of lowest and of highest degree of the squared polynomials would give nonzero terms of degree other than \(2d\).
We write \(\DSOSForms\) and \(\SDSOSForms\) for the sets of DSOS and of SDSOS forms, and \(\SOSForms\) for the set of sums of squares of polynomials.
A symmetric matrix \(\Gamma\) indexed by \(\Naturals^{n}_{d}\) is a Gram matrix of \(q\) if \(q=\left[\FormVariable\right]_{d}^{\Transpose}\Gamma\left[\FormVariable\right]_{d}\).
The squares \(\left(\FormVariable^{\nu}\right)^{2}\), \(\left(\FormVariable^{\nu}\pm \FormVariable^{\nu^{\prime}}\right)^{2}\) and \(\left(a \FormVariable^{\nu}+b\FormVariable^{\nu^{\prime}}\right)^{2}\) have the Gram matrices \(e_{\nu} e_{\nu}^{\Transpose}\), \(\left(e_{\nu}\pm e_{\nu^{\prime}}\right)\left(e_{\nu}\pm e_{\nu^{\prime}}\right)^{\Transpose}\) and \(vv^{\Transpose}\) with \(v=a e_{\nu}+be_{\nu^{\prime}}\).
The first two are diagonally dominant, and the display above writes every diagonally dominant matrix as a nonnegative combination of them; a positive semidefinite matrix that vanishes outside the rows and columns of \(\nu\) and \(\nu^{\prime}\) is, by its spectral decomposition, a sum of matrices \(vv^{\Transpose}\) of the third kind.
Hence \(q\in\DSOSForms\) if and only if \(q\) has a diagonally dominant Gram matrix, and \(q\in\SDSOSForms\) if and only if \(q\) has a Gram matrix that is a sum of \(2\times2\) blocks.
Likewise, if \(q\) has a Gram matrix \(\Gamma=\sum_{k}v_{k}v_{k}^{\Transpose}\succeq0\), then \(q=\sum_{k}\left(v_{k}^{\Transpose}\left[\FormVariable\right]_{d}\right)^{2}\) lies in \(\SOSForms\).
We put
\[
\begin{aligned}
\PolyaDSOSCone[n]{r}&\coloneqq\left\{M\in\Sym[n]:\left(\sum_{i}\FormVariable_{i}^{2}\right)^{r}p_{M}\in\DSOSForms\right\},\\
\PolyaSDSOSCone[n]{r}&\coloneqq\left\{M\in\Sym[n]:\left(\sum_{i}\FormVariable_{i}^{2}\right)^{r}p_{M}\in\SDSOSForms\right\},\\
\PolyaSOSCone[n]{r}&\coloneqq\left\{M\in\Sym[n]:\left(\sum_{i}\FormVariable_{i}^{2}\right)^{r}p_{M}\in\SOSForms\right\}.
\end{aligned}
\]
The last is the cone of \citep{parrilo2000structured} in the form \citep[Definition 3.1]{deklerk2002approximation}, and \(\PolyaSOSCone[n]{1}\) is the cone of \hyperref[sec-rho]{Section~\ref*{sec-rho}}.
Let \(\PolyaDSOS{r}\left(G\right)\), \(\PolyaSDSOS{r}\left(G\right)\) and \(\PolyaSOS{r}\left(G\right)\) be the optimal values of \hyperref[eq-zeta]{(\ref*{eq-zeta})} with \(\PolyaDSOSCone[n]{r}\), \(\PolyaSDSOSCone[n]{r}\) and \(\PolyaSOSCone[n]{r}\) in place of \(\PolyaCone[n]{r}\), with the value \(\infty\) when the program is infeasible.
Diagonal dominance of \(\Gamma\) is a system of linear inequalities in \(\Gamma\) and in variables \(t_{ij}\ge|\Gamma_{ij}|\); thus \(\PolyaDSOS{r}\left(G\right)\) is the value of a linear program.

\begin{lemma}[]\protect\hypertarget{lem-gram}{}\label{lem-gram}

For \(M\in\Sym[n]\), \(M\in\PolyaDSOSCone[n]{0}\) if and only if \(M=\Phi+Z\) with \(\Phi\) diagonally dominant and \(Z\in\NN[n]\), and \(M\in\PolyaSDSOSCone[n]{0}\) if and only if \(M=\Phi+Z\) with \(\Phi\) a sum of \(2\times2\) blocks and \(Z\in\NN[n]\).

\end{lemma}

\begin{proof}
The monomials of degree \(2\) are the \(\FormVariable_{i}^{2}\), indexed by \(2e_{i}\), and the \(\FormVariable_{i}\FormVariable_{j}\) with \(i<j\), indexed by \(e_{i}+e_{j}\).
Let \(\Gamma\) be a Gram matrix of \(p_{M}\).
A product of two monomials of degree \(2\) equals \(\FormVariable_{i}^{4}\) only for \(\FormVariable_{i}^{2}\cdot \FormVariable_{i}^{2}\), and equals \(\FormVariable_{i}^{2}\FormVariable_{j}^{2}\), with \(i\ne j\), only for \(\FormVariable_{i}^{2}\cdot \FormVariable_{j}^{2}\), \(\FormVariable_{j}^{2}\cdot \FormVariable_{i}^{2}\) and \(\FormVariable_{i}\FormVariable_{j}\cdot \FormVariable_{i}\FormVariable_{j}\).
Comparing the coefficients of these monomials gives
\[
M_{ii}=\Gamma_{2e_{i},2e_{i}},\qquad 2M_{ij}=2\Gamma_{2e_{i},2e_{j}}+\Gamma_{e_{i}+e_{j},e_{i}+e_{j}}\quad\left(i\ne j\right).
\]
Hence \(M=\Gamma^{\prime}+Z\), where \(\Gamma^{\prime}\) is the principal submatrix of \(\Gamma\) on the indices \(2e_{i}\), relabeled by \(i\), and \(Z\) has zero diagonal and \(Z_{ij}=\frac12\Gamma_{e_{i}+e_{j},e_{i}+e_{j}}\) for \(i\ne j\).
If \(\Gamma\) is diagonally dominant, or a sum of \(2\times2\) blocks, then so is \(\Gamma^{\prime}\), and \(Z\ge0\) because \(\Gamma\succeq0\) has a nonnegative diagonal.

Conversely, let \(M=\Phi+Z\) with \(Z\in\NN[n]\) and \(\Phi\) diagonally dominant, or a sum of \(2\times2\) blocks.
Then
\[
p_{M}\left(\FormVariable\right)=\sum_{i,j}\Phi_{ij}\FormVariable_{i}^{2}\FormVariable_{j}^{2}+\sum_{i}Z_{ii}\FormVariable_{i}^{4}+\sum_{i<j}2Z_{ij}\left(\FormVariable_{i}\FormVariable_{j}\right)^{2},
\]
and thus \(p_{M}=\left[\FormVariable\right]_{2}^{\Transpose}\Gamma\left[\FormVariable\right]_{2}\) for the matrix \(\Gamma\) that equals \(\Phi+\Diag\left(Z_{11},\dots,Z_{nn}\right)\) on the indices \(2e_{i}\), has the diagonal entries \(2Z_{ij}\ge0\) at the indices \(e_{i}+e_{j}\), and vanishes elsewhere.
This matrix has the same property as \(\Phi\).
\end{proof}

\begin{proposition}[]\protect\hypertarget{prp-dsos-chain}{}\label{prp-dsos-chain}

For every \(r\ge0\),
\[
\PolyaCone[n]{r}\subseteq\PolyaDSOSCone[n]{r}\subseteq\PolyaSDSOSCone[n]{r}\subseteq\PolyaSOSCone[n]{r}\subseteq\COP[n],\qquad\PolyaDSOSCone[n]{r}\subseteq\PolyaDSOSCone[n]{r+1},\qquad\PolyaSDSOSCone[n]{r}\subseteq\PolyaSDSOSCone[n]{r+1}.
\]
Hence \(\tfrac12\sqrt{\Biclique\left(G\right)}\le\PolyaSOS{r}\left(G\right)\le\PolyaSDSOS{r}\left(G\right)\le\PolyaDSOS{r}\left(G\right)\le\Polya{r}\left(G\right)\), the bounds \(\PolyaDSOS{r}\left(G\right)\) and \(\PolyaSDSOS{r}\left(G\right)\) are nonincreasing in \(r\), and if \(\Biclique\left(G\right)\ge1\) they converge to \(\tfrac12\sqrt{\Biclique\left(G\right)}\).

\end{proposition}

\begin{proof}
Let \(M\in\PolyaCone[n]{r}\).
By \hyperref[sec-polya]{Section~\ref*{sec-polya}}, \(\left(\sum_{i}\FormVariable_{i}^{2}\right)^{r}p_{M}\) has only the monomials \(\FormVariable^{2\nu}\), \(\nu\in\Naturals^{n}_{r+2}\), and their coefficients \(c_{\nu}\) are nonnegative; hence
\[
\left(\sum_{i}\FormVariable_{i}^{2}\right)^{r}p_{M}\left(\FormVariable\right)=\sum_{\nu} c_{\nu}\left(\FormVariable^{\nu}\right)^{2}=\left[\FormVariable\right]_{r+2}^{\Transpose}\Diag\left(c\right)\left[\FormVariable\right]_{r+2},
\]
and \(\Diag\left(c\right)\) is diagonally dominant; thus \(M\in\PolyaDSOSCone[n]{r}\).
This inclusion is observed in \citep{ahmadi2019dsos}.
The next two inclusions follow from the inclusions between diagonally dominant matrices, sums of \(2\times2\) blocks and positive semidefinite matrices, and from the remark on sums of squares above.
If \(M\in\PolyaSOSCone[n]{r}\), then \(\left(\sum_{i}\FormVariable_{i}^{2}\right)^{r}p_{M}\left(\FormVariable\right)\ge0\) for every \(\FormVariable\); hence \(p_{M}\ge0\) and \(M\in\COP[n]\).

For the nesting, let \(q=\left(\sum_{i}\FormVariable_{i}^{2}\right)^{r}p_{M}=\left[\FormVariable\right]_{d}^{\Transpose}\Gamma\left[\FormVariable\right]_{d}\) with \(d=r+2\) and \(\Gamma\) diagonally dominant, or a sum of \(2\times2\) blocks.
For each \(k\), the entries of \(\FormVariable_{k}\left[\FormVariable\right]_{d}\) are the distinct monomials \(\FormVariable^{\nu+e_{k}}\), \(\nu\in\Naturals^{n}_{d}\); thus \(\FormVariable_{k}^{2}q=\left[\FormVariable\right]_{d+1}^{\Transpose}\Gamma_{k}\left[\FormVariable\right]_{d+1}\), where \(\Gamma_{k}\) has the entry \(\Gamma_{\nu\nu^{\prime}}\) at the indices \(\nu+e_{k}\), \(\nu^{\prime}+e_{k}\) and vanishes elsewhere.
Each \(\Gamma_{k}\) is \(\Gamma\) relabeled and bordered by zeros; hence \(\sum_{k}\Gamma_{k}\) has the same property as \(\Gamma\), and it is a Gram matrix of \(\left(\sum_{k}\FormVariable_{k}^{2}\right)q=\left(\sum_{i}\FormVariable_{i}^{2}\right)^{r+1}p_{M}\).

A larger cone gives a smaller value in \hyperref[eq-zeta]{(\ref*{eq-zeta})}; thus the inclusions give all the inequalities but the first, and the nesting gives the monotonicity in \(r\).
For the first, a feasible point for \(\PolyaSOS{r}\) is feasible for \hyperref[eq-dual]{(\ref*{eq-dual})}, whose value is \(\tfrac12\sqrt \Biclique\) by \hyperref[cor-dual]{Corollary~\ref*{cor-dual}}.
If \(\Biclique\ge1\), then \(\Polya{r}\) decreases to \(\tfrac12\sqrt \Biclique\) by \hyperref[prp-polya-converge]{Proposition~\ref*{prp-polya-converge}}, and so do \(\PolyaDSOS{r}\) and \(\PolyaSDSOS{r}\), which are nonincreasing and lie between \(\tfrac12\sqrt \Biclique\) and \(\Polya{r}\).
\end{proof}

At level \(0\) the DSOS bound has a closed form.
Let \(\hat d\) be the maximum degree of \(\Complement{G}\).
For \(M=\lambda \Identity-\CostMatrix+\sum_{ij\in E}\EdgeMultiplier_{ij}S_{ij}\), we have \(M_{ii}=\lambda\), \(M_{ij}=-\frac12\) for the edges \(ij\) of \(\Complement{G}\), \(M_{ij}=\EdgeMultiplier_{ij}-\frac12\) for \(ij\in E\), and \(M_{ij}=0\) for distinct \(i\), \(j\) on the same side.

\begin{proposition}[]\protect\hypertarget{prp-dsos-degree}{}\label{prp-dsos-degree}

For every bipartite graph \(G\), \(\PolyaDSOS{0}\left(G\right)=\frac12\hat d\), attained by \(\lambda=\frac12\hat d\) and \(\EdgeMultiplier_{ij}=\frac12\) for all \(ij\in E\).
Moreover \(\Biclique\left(G\right)\le\hat d^{2}\); hence \(\PolyaDSOS{0}\left(G\right)=\tfrac12\sqrt{\Biclique\left(G\right)}\) if and only if \(\hat d=0\) or \(G\) has a biindependent pair \(\left(A,B\right)\) with \(|A|=|B|=\hat d\).

\end{proposition}

\begin{proof}
Let \(M=\lambda \Identity-\CostMatrix+\sum_{ij\in E}\EdgeMultiplier_{ij}S_{ij}\in\PolyaDSOSCone[n]{0}\).
By \hyperref[lem-gram]{Lemma~\ref*{lem-gram}}, \(M=\Phi+Z\) with \(\Phi\) diagonally dominant and \(Z\in\NN[n]\); thus \(\Phi_{ii}\le\lambda\), and \(\Phi_{ij}\le-\frac12\) for every edge \(ij\) of \(\Complement{G}\).
For a vertex \(i\) with \(\hat d\) neighbors in \(\Complement{G}\), row \(i\) of \(\Phi\) has \(\hat d\) off-diagonal entries of absolute value at least \(\frac12\), and diagonal dominance gives \(\lambda\ge\Phi_{ii}\ge\frac12\hat d\).
Conversely, for \(\lambda=\frac12\hat d\) and \(\EdgeMultiplier_{ij}=\frac12\), the off-diagonal entries of \(M\) are \(-\frac12\) on the edges of \(\Complement{G}\) and \(0\) elsewhere; hence the off-diagonal entries of row \(i\) have absolute values summing to half the degree of \(i\) in \(\Complement{G}\), which is at most \(\lambda=M_{ii}\).
Thus \(M\) is diagonally dominant, and \(M\in\PolyaDSOSCone[n]{0}\) by \hyperref[lem-gram]{Lemma~\ref*{lem-gram}}.

Let \(\left(A,B\right)\) be biindependent with \(A\) and \(B\) nonempty.
A vertex of \(A\) is adjacent in \(\Complement{G}\) to every vertex of \(B\); thus \(|B|\le\hat d\), and likewise \(|A|\le\hat d\).
Hence \(\Biclique\le\hat d^{2}\), which also holds if \(\Biclique=0\), and if \(\hat d\ge1\), then \(\Biclique=\hat d^{2}\) if and only if some biindependent pair has \(|A|=|B|=\hat d\).
If \(\hat d=0\), then \(\Complement{G}\) has no edges; thus \(G\) is complete bipartite and \(\Biclique=0=\hat d^{2}\).
As \(\PolyaDSOS{0}=\frac12\hat d\), it equals \(\tfrac12\sqrt \Biclique\) if and only if \(\Biclique=\hat d^{2}\).
\end{proof}

For example, if \(G\) is the crown \(K_{s,s}\) minus a perfect matching, \(s\ge2\), then \(\Complement{G}\) is a perfect matching; hence \(\hat d=\Biclique=1\) and \(\PolyaDSOS{0}\left(G\right)=\frac12=\tfrac12\sqrt \Biclique\), while the Pólya bounds are never exact (\hyperref[thm-polya-inexact]{Theorem~\ref*{thm-polya-inexact}}).

\begin{proposition}[]\protect\hypertarget{prp-dsos-spectral}{}\label{prp-dsos-spectral}

For every bipartite graph \(G\),
\[
\tfrac12\|\OnesMatrix-\Biadjacency\|_{2}\le\PolyaSDSOS{0}\left(G\right)\le\PolyaDSOS{0}\left(G\right)=\tfrac12\hat d .
\]
Thus the penalty \(\gamma=\frac12+\PolyaSDSOS{0}\left(G\right)\), which satisfies the hypothesis of \hyperref[thm-penalty]{Theorem~\ref*{thm-penalty}} by \hyperref[prp-dsos-chain]{Proposition~\ref*{prp-dsos-chain}}, is never smaller than the spectral threshold \(\frac12\left(1+\|\OnesMatrix-\Biadjacency\|_{2}\right)\).

\end{proposition}

\begin{proof}
The second inequality and the equality are \hyperref[prp-dsos-chain]{Proposition~\ref*{prp-dsos-chain}} and \hyperref[prp-dsos-degree]{Proposition~\ref*{prp-dsos-degree}}.
Let \(M=\lambda \Identity-\CostMatrix+\sum_{ij\in E}\EdgeMultiplier_{ij}S_{ij}\in\PolyaSDSOSCone[n]{0}\).
By \hyperref[lem-gram]{Lemma~\ref*{lem-gram}}, \(M=\Phi+Z\) with \(\Phi\) a sum of \(2\times2\) blocks and \(Z\in\NN[n]\).
Grouping the terms of \(\Phi\) by the indices outside which they vanish, and assigning a term that vanishes outside a single index \(i\) to some pair containing \(i\), we get \(\Phi=\sum_{i<j}\Phi^{ij}\) with \(\Phi^{ij}\succeq0\) zero outside the rows and columns \(i\) and \(j\); then \(\Phi_{ij}=\Phi^{ij}_{ij}\) for \(i\ne j\).
Let \(\Psi^{ij}\) be \(\Phi^{ij}\) with its two off-diagonal entries replaced by \(\min\left\{\Phi_{ij},0\right\}\).
It equals \(\Phi^{ij}\) if \(\Phi_{ij}\le0\) and is diagonal with nonnegative entries otherwise; thus \(\Psi^{ij}\succeq0\).
Hence \(\Psi=\sum_{i<j}\Psi^{ij}\) is a sum of \(2\times2\) blocks, with \(\Psi_{ii}=\Phi_{ii}\) and \(\Psi_{ij}=\min\left\{\Phi_{ij},0\right\}\) for \(i\ne j\), and \(M=\Psi+Z^{\prime}\) with \(Z^{\prime}=Z+\Phi-\Psi\in\NN[n]\).
Therefore \(\Psi_{ii}\le M_{ii}=\lambda\), \(\Psi_{ij}\le0\) for \(i\ne j\), and \(\Psi_{ij}\le M_{ij}=-\frac12\) for every edge \(ij\) of \(\Complement{G}\).

The matrix \(2\CostMatrix-\Adjacency\) is the adjacency matrix of \(\Complement{G}\), and its largest eigenvalue is \(\|\OnesMatrix-\Biadjacency\|_{2}\), as noted after \hyperref[thm-penalty]{Theorem~\ref*{thm-penalty}}.
Let \(u\) be a unit eigenvector for this eigenvalue, and \(\left|u\right|\) the vector of the absolute values of its entries.
As \(2\CostMatrix-\Adjacency\ge0\), \(\left|u\right|\) is a nonnegative unit vector with \(\left|u\right|^{\Transpose}\left(2\CostMatrix-\Adjacency\right)\left|u\right|\ge u^{\Transpose}\left(2\CostMatrix-\Adjacency\right)u=\|\OnesMatrix-\Biadjacency\|_{2}\).
Since \(\Psi\succeq0\),
\[
\begin{aligned}
0\le\left|u\right|^{\Transpose}\Psi\left|u\right|&=\sum_{i}\Psi_{ii}u_{i}^{2}+\sum_{i\ne j}\Psi_{ij}\left|u_{i}\right|\left|u_{j}\right|\\
&\le\lambda-\tfrac12\left|u\right|^{\Transpose}\left(2\CostMatrix-\Adjacency\right)\left|u\right|\le\lambda-\tfrac12\|\OnesMatrix-\Biadjacency\|_{2} .
\end{aligned}
\]
Thus every feasible \(\lambda\) is at least \(\frac12\|\OnesMatrix-\Biadjacency\|_{2}\).
\end{proof}

\hyperref[tbl-dsos]{Table~\ref*{tbl-dsos}} gives \(\PolyaDSOS{r}\) for \(r\le2\) on the graphs of \hyperref[fig-polya]{Figure~\ref*{fig-polya}} and on \(N_{2}\).
We solved the linear programs as in \hyperref[sec-numerics]{Section~\ref*{sec-numerics}}, with \(\Gamma\) indexed by \(\Naturals^{n}_{r+2}\); the entries \(\Gamma_{\nu\nu^{\prime}}\) with \(\nu+\nu^{\prime}\notin2\Naturals^{n}\) multiply only monomials that do not occur in \(\left(\sum_{i}\FormVariable_{i}^{2}\right)^{r}p_{M}\), and setting them to \(0\) keeps \(\Gamma\) diagonally dominant; hence we omitted them.
On \(3K_{2}\) the Pólya program is infeasible for \(r\le1\) and \(\Polya{2}\left(3K_{2}\right)=2=2\sqrt2\cdot\tfrac12\sqrt \Biclique\) (\hyperref[fig-polya]{Figure~\ref*{fig-polya}}), while \(\PolyaDSOS{0}\left(3K_{2}\right)=1\) already.
Unlike the Pólya bounds, the DSOS bounds also distinguish \(3K_{2}\) from its induced subgraph with one edge on \(2+2\) vertices.

\begin{table}[H]
\centering\small
\caption{The DSOS bounds $\PolyaDSOS{r}$ on the graphs of {\hyperref[fig-polya]{Figure~\ref*{fig-polya}}} and on the graph $N_{2}$ of {\hyperref[prp-sharp]{Proposition~\ref*{prp-sharp}}}.}\label{tbl-dsos}
\newcommand{\dsosrow}[3]{\readrows#1{filter strcmp={\colgraph}{#2}}{\TableData{dsos}}{#3 & $\fixed{4}{\colhalfsqrtg}$ & $\fixed{4}{\collevelzero}$ & $\fixed{4}{\collevelone}$ & $\fixed{4}{\colleveltwo}$}}
\dsosrow\rowmthree{matching-3}{$3K_{2}$}
\dsosrow\rowoneedge{lpv-example-1}{one edge on $2+2$ vertices}
\dsosrow\rowmfour{matching-4}{$4K_{2}$}
\dsosrow\rowcrown{crown-4}{$K_{4,4}$ minus a perfect matching}
\dsosrow\rownk{near-biclique-2}{$N_{2}$}

\begin{tabular}{@{}lrrrr@{}}
\toprule
& & \multicolumn{3}{c}{$\PolyaDSOS{r}$}\\
\cmidrule(l){3-5}
$G$ & $\frac12\sqrt \Biclique$ & $r=0$ & $r=1$ & $r=2$\\
\midrule
\rowmthree\rowoneedge\rowmfour\rowcrown\rownk
\bottomrule
\end{tabular}
\end{table}

\section{Computational Details}\label{sec-comp}

\textbf{The cone \(\PolyaSOSCone[n]{1}\).}
By \citep{parrilo2000structured} in the form \citep[(13)--(16)]{deklerk2002approximation}, a matrix \(M\in\Sym[n]\) lies in \(\PolyaSOSCone[n]{1}\) if and only if there are matrices \(M^{\left(1\right)},\dots,M^{\left(n\right)}\in\Sym[n]\) such that
\[
\begin{aligned}
&M-M^{\left(i\right)}\succeq0, &\quad& M^{\left(i\right)}_{ii}=0, && i=1,\dots,n,\\
&M^{\left(j\right)}_{ii}+2M^{\left(i\right)}_{ij}=0, && && i\ne j,\\
&M^{\left(i\right)}_{jk}+M^{\left(j\right)}_{ik}+M^{\left(k\right)}_{ij}\ge0, && && i<j<k .
\end{aligned}
\]
Thus \(\PolyaSOS{1}\) is a semidefinite program with \(n\) blocks of size \(n\) and about \(\frac16n^{3}\) linear inequalities.

\begin{table}[tbp]
\centering\footnotesize\setlength{\tabcolsep}{3.2pt}
\caption{The bounds $\LasserreBicliqueRatio{1}$, $\ConeBicliqueRatio{\DNN}$ and $\PolyaSOS{1}$ on the random graphs of {\hyperref[tbl-random]{Table~\ref*{tbl-random}}}, for each size and each edge probability $p$; \# is the number of graphs.
Columns 4--6: number of graphs on which the bound is within $\sci{\Value{exacttol}}$ of $\frac12\sqrt \Biclique$.
Columns 7--12: gap of the bound relative to $\frac12\sqrt \Biclique$, in percent, mean and maximum over the graphs of the row.}\label{tbl-random-p}
\readrows\balancedrows{filter strcmp={\colblock}{balanced}, and filter not strcmp={\colp}{all}}{\TableData{bounds}}{\ifnum\thecsvrow>1 \ifcsvstrcmp{\colp}{0.3}{\csvexpnot\addlinespace}{}\fi\ifcsvstrcmp{\colp}{0.3}{$\colsize$}{} & $\colp$ & \boundscells}
\readrows\unbalancedrows{filter strcmp={\colblock}{unbalanced}, and filter not strcmp={\colp}{all}}{\TableData{bounds}}{\ifnum\thecsvrow>1 \ifcsvstrcmp{\colp}{0.3}{\csvexpnot\addlinespace}{}\fi\ifcsvstrcmp{\colp}{0.3}{$\colsize$}{} & $\colp$ & \boundscells}

\begin{tabular}{@{}rrrrrrrrrrrr@{}}
\toprule
& & & \multicolumn{3}{c}{exact} & \multicolumn{3}{c}{mean gap (\%)} & \multicolumn{3}{c}{max gap (\%)}\\
\cmidrule(lr){4-6}\cmidrule(lr){7-9}\cmidrule(l){10-12}
$m_{1}+m_{2}$ & $p$ & \boundsheads\\
\midrule
\multicolumn{12}{@{}l}{\emph{balanced, $m_{1}=m_{2}$}}\\
\balancedrows
\midrule
\multicolumn{12}{@{}l}{\emph{unbalanced, $m_{1}<m_{2}$}}\\
\unbalancedrows
\bottomrule
\end{tabular}
\end{table}

\textbf{Exact evaluation of \(\QuadraticMax{\gamma}\).}
Let \(Q=\CostMatrix-\gamma \Adjacency\), and let \(u\) be a maximizer of \(u^{\Transpose}Qu\) over nonnegative unit vectors, with support \(S\).
Then \(u_{S}\) is a positive local maximizer of \(v^{\Transpose}Q_{S}v\) on the unit sphere of \(\Reals^{S}\); thus \(Q_{S}u_{S}=\mu u_{S}\) with \(\mu=u^{\Transpose}Qu\), and \(\mu=\lambda_{\max}\left(Q_{S}\right)\), because a local maximizer of a quadratic form on the unit sphere is a global one.
Conversely, a positive eigenvector of a principal submatrix \(Q_{S}\), extended by zeros, is feasible with value its eigenvalue.
Hence \(\QuadraticMax{\gamma}\left(G\right)\) is the maximum, over the supports \(S\), of the eigenvalues of \(Q_{S}\) whose eigenspace contains a positive vector.
We find it by enumerating the \(2^{n}-1\) supports, and decide whether an eigenspace contains a positive vector by a small linear program, which is exact also for repeated eigenvalues.
Since \(\gamma\mapsto\QuadraticMax{\gamma}\) is nonincreasing by \hyperref[eq-mu-sigma]{(\ref*{eq-mu-sigma})} and \hyperref[lem-sigma]{Lemma~\ref*{lem-sigma}}, \(\SmallestPenalty\left(G\right)\) is found by bisection, stopped at an interval of length \(\sci{\Value{bisectiontol}}\).
The bisection counts a penalty \(\gamma\) as exact when the computed \(\QuadraticMax{\gamma}\) exceeds \(\tfrac12\sqrt \Biclique\) by at most \(\sci{\Value{penaltytol}}\).
The computed value can therefore lie below \(\SmallestPenalty\) by more than the length of this interval when \(\QuadraticMax{\gamma}\) approaches \(\tfrac12\sqrt \Biclique\) with zero slope, as it does on the graphs \(N_{k}\).

\textbf{Random graphs and tolerances.}
A random bipartite graph with sides \(m_{1}\le m_{2}\) and edge probability \(p\) has the biadjacency matrix \texttt{rng.random((m1,\ m2))\ \textless{}\ p}, where \texttt{rng} is \texttt{numpy.random.default\_rng} seeded with \(1000m_{1}+100\lfloor10p\rfloor+s\) if \(m_{1}=m_{2}\) and with \(100000m_{2}+1000m_{1}+100\lfloor10p\rfloor+s\) otherwise, \(s=0,\dots,99\).
For the penalty threshold, every size uses the second seed, and the graphs without edges and the complete ones are skipped.
\texttt{MOSEK} 11.2 was run with relative gap and feasibility tolerances \(\sci{\Value{solvertol}}\).
On each of the \(\Value{rhooneinexact}\) graphs on which \(\PolyaSOS{1}\) is not exact, it reports an optimal solution whose primal and dual objective values agree to within \(\sciatmost{1}{\Value{dualgap}}\).

\textbf{Penalties below the threshold.}
\hyperref[tbl-ascent-below]{Table~\ref*{tbl-ascent-below}} runs the local ascent of \hyperref[sec-ascent]{Section~\ref*{sec-ascent}} from the random starts with the penalties \(\frac12\left(1+c\sqrt \Biclique\right)\) for \(c=\frac14\), \(\frac12\) and \(\frac34\), which lie below the threshold \(\PenaltyThreshold{g}\) of \hyperref[sec-exp-threshold]{Section~\ref*{sec-exp-threshold}}.
With \(c=\frac34\) the ascent succeeds more often than with \(\PenaltyThreshold{g}\) on seven of the nine size and density classes.
With \(c=\frac12\), between \(\fixed{0}{\Value{halfedgemin}}\%\) and \(\fixed{0}{\Value{halfedgemax}}\%\) of the runs at \(p\le0.5\) end with an edge in the support, and with \(c=\frac14\) at least \(\atleast{0}{\Value{quarteredgemin}}\%\) of all runs do, and hence success then depends on the repair; after it, \(c=\frac12\) does best on five of the nine classes, and no single value of \(c\) is best throughout.
The smallest exact penalty \(\SmallestPenalty\) cannot be computed by enumeration at these sizes; hence we could not test it.

\begin{table}[tbp]
\centering\small\setlength{\tabcolsep}{4.5pt}
\caption{Fraction of random starts from which the local ascent ends at a maximum pair, without the greedy repair (ends) and with it (repaired), for the penalties $\gamma=\frac12\left(1+c\sqrt \Biclique\right)$ on random graphs of size $G$ and edge probability $p$; $c=1$ gives the threshold of {\hyperref[thm-penalty]{Theorem~\ref*{thm-penalty}}}.}\label{tbl-ascent-below}
\readrows\randomrows{filter strcmp={\colstarts}{random}}{\TableData{ascent}}{$\colsize$ & $\colp$ & $\fixed{3}{\colcquartersuccess}$ & $\fixed{3}{\colcquarterrepaired}$ & $\fixed{3}{\colchalfsuccess}$ & $\fixed{3}{\colchalfrepaired}$ & $\fixed{3}{\colcthreequarterssuccess}$ & $\fixed{3}{\colcthreequartersrepaired}$ & $\fixed{3}{\colgammaonesuccess}$ & $\fixed{3}{\colgammaonerepaired}$}

\begin{tabular}{@{}lrrrrrrrrr@{}}
\toprule
& & \multicolumn{2}{c}{$c=\frac14$} & \multicolumn{2}{c}{$c=\frac12$} & \multicolumn{2}{c}{$c=\frac34$} & \multicolumn{2}{c}{$c=1$}\\
\cmidrule(lr){3-4}\cmidrule(lr){5-6}\cmidrule(lr){7-8}\cmidrule(l){9-10}
$G$ & $p$ & ends & repaired & ends & repaired & ends & repaired & ends & repaired\\
\midrule
\randomrows
\bottomrule
\end{tabular}
\end{table}

\renewcommand\refname{References}
\bibliography{refs}

\end{document}